\documentclass[12pt]{article}
\usepackage{makeidx}
\usepackage{amssymb,amsthm,amsmath,amsfonts,enumerate}
\usepackage{setspace,graphicx,caption,subcaption,float,relsize}
\usepackage{url}
\usepackage{hyperref}
\usepackage{rotating}
\usepackage{natbib}
\usepackage{array}

\hypersetup{
  colorlinks = true, urlcolor = cyan, linkcolor = blue, citecolor = red }
\renewenvironment{abstract}
 {\small
  \begin{center}
  \bfseries \abstractname\vspace{-.0em}\vspace{0pt}
  \end{center}
  \list{}{    \setlength{\leftmargin}{0mm}
    \setlength{\rightmargin}{\leftmargin}  }  \item\relax}
 {\endlist}
\makeatletter
\def\maketag@@@#1{\hbox{\m@th\normalfont\normalsize#1}}
\makeatother

\allowdisplaybreaks
\newtheorem {theorem}{Theorem}[section]
\newtheorem {assumption}{Assumption}

\newtheorem{lemma}[theorem]{Lemma}

\newtheorem{remark}{Remark}[section]

\allowdisplaybreaks
\input{tcilatex}
\begin{document}

\title{Testing Many Zero Restrictions in a High Dimensional Linear
Regression Setting}
\author{Jonathan B. Hill\thanks{%
Department of Economics, University of North Carolina, Chapel Hill, North
Carolina, U.S. E-mail:\texttt{jbhill@email.unc.edu}; \texttt{%
https://jbhill.web.unc.edu}. This paper was previously circulated under the
title \textquotedblleft \textit{Testing (Infinitely) Many Zero Restrictions}%
".}\medskip \\
Dept. of Economics, University of North Carolina, U.S.\medskip \\
\large{Original draft: January 2023}}
\date{{\large This draft:} \today}
\maketitle

\begin{abstract}
We propose a test of many zero parameter restrictions in a high dimensional
linear iid regression model with $k$ $>>$ $n$ regressors. The test statistic
is formed by estimating key parameters one at a time based on many low
dimension regression models with nuisance terms. The parsimoniously
parametrized models identify whether the original parameter of interest is
or is not zero. Estimating fixed low dimension sub-parameters ensures
greater estimator accuracy, it does not require a sparsity assumption nor
therefore a regularized estimator, it is computationally fast compared to,
e.g., de-biased Lasso, and using only the largest in a sequence of weighted
estimators reduces test statistic complexity and therefore estimation error.
We provide a parametric wild bootstrap for p-value computation, and prove
the test is consistent and has non-trivial $\sqrt{n/\{\ln (n)\mathcal{M}%
_{n}\}}$-local-to-null power where $\mathcal{M}_{n}$ is the $l_{\infty }$
covariate fourth moment.\bigskip \newline
\end{abstract}

\doublespacing

\setstretch{1.7}

\section{Introduction\label{sec:intro}}

Regression settings where the number $k$ of covariates $x_{t}$\ may be much
larger than the sample size $n$ ($k>>n$) is a natural possibility in cross
sections, panels and spatial settings where an enormous amount of
information is available. In economics and finance, for example, reasons
include the size of surveys (e.g. U.S. Census, Current Population Survey,
National Longitudinal Survey of Youth), increasingly sophisticated survey
techniques (e.g. text records and word counts, household transaction data),
the unit of observation (e.g. county wide house prices are dependent on
neighbor county prices), use of dummy variables for groups (e.g. age, race,
region, etc.), and panel data fixed effects. See, for instance, \cite%
{FanLvQi2011} and \cite{BelloniChernozhukovHansen2014} for recent surveys.
In non-social sciences high dimensionality arises with, for instance,
genomic data and network and communication data.\footnote{%
As a testament to the increasing amount of available data, in November 2022
the 27th General Conference on Weights and Measures presented new numerical
prefixes handling, among other magnitudes, $10^{30}$.}

Inference in this case is typically done post estimation based on a
regression or moment condition model with an imposed sparsity condition, or
related structural assumption. Sparsity in practice is achieved by a
shrinkage or regularized estimator, while recent theory allows for an
increasing number of candidate and relevant covariates as the sample size
increases \citep[see, e.g.,][]{FanPeng2004,HuangHorowitzMa2008}. Valid
inference, however, typically exists only for the non-zero valued parameter
subset, although there is a nascent literature for gaining inference on the
(number of) zero valued parameters \citep[e.g.][]{Carpentier_Verzelen_2019}.
Conversely, \cite{Buhlmann2013} develops a test of zero restrictions for a
high dimensional linear model with fixed design $\{x_{1},...,x_{n}\}^{\prime
}$ using bias-corrected ridge regression. Their p-value construction,
however, is based on a linear model with Gaussian errors and an upper bound
on the distribution of the ridge estimator. Low dimension projection methods
offer another solution \citep[e.g.][]{ZhangZhang2014}. Imposing sparsity at
all, however, may be too restrictive, in particular without justification
via a pre-estimation test procedure.

\cite{Dezeure_etal_2017} work in the de-biased Lasso [DBL] setting of \cite%
{ZhangZhang2014} to deliver a useful post-estimation max-t-test under a
shrinkage assumption. The max-test in \cite{Dezeure_etal_2017} is similar to
ours in the sense that they bootstrap the maximum t-statistic over an
arbitrary $k$-set of t-statistics, as long as $\ln (k)$ $=$ $o(n^{1/7})$
under covariate boundedness. Boundedness, however, is unlikely to exist in
key settings in economics, including nonlinear models involving propensity
scores, inverse probability weighting, and log-linearization; and many
standard models in finance like ARMA-GARCH. Although we work in an iid
linear setting to focus ideas, our methods can only reasonably, and
usefully, be extended to time series models under unboundedness.

Consider a sample of iid observations $\{x_{t},y_{t}\}_{t=1}^{n}$ with
scalar $y_{t}$, and covariates or design points $x_{t}$ $=$ $[x_{\delta
,t}^{\prime },x_{\theta ,t}^{\prime }]^{\prime }$, and a linear regression
model for simplicity and to fix ideas:%
\begin{equation}
y_{t}=\delta _{0}^{\prime }x_{\delta ,t}+\theta _{0}^{\prime }x_{\theta
,t}+\epsilon _{t}=\beta _{0}^{\prime }x_{t}+\epsilon _{t}\text{.}
\label{model}
\end{equation}%
This paper presents a test of $H_{0}$ $:$ $\theta _{0}$ $=$ $0$ vs. $H_{1}$ $%
:$ $\theta _{0,i}$ $\neq $ $0$ for at least one $i$. The nuisance parameter $%
\delta _{0}$ has a fixed dimension $k_{\delta }$ $<$ $n$, but $\theta _{0}$
may have dimension $k_{\theta }$ $>>$ $n$. We assume $k_{\theta }$ $%
\rightarrow $ $\infty $, but $k_{\theta }$ $\rightarrow $ $k_{\theta }^{\ast
}$ $\in $ $\mathbb{N}$ can be easily allowed. In order to show that $%
k_{\theta }$ depends on $n$ we write $k_{\theta ,n}$ $=$ $k_{\theta }$
throughout. Let $\beta _{0}$ be the unique minimizer of the squared error
loss $E[(y_{t}$ $-$ $\beta ^{\prime }x_{t})^{2}]$, where $E[\epsilon _{t}]$ $%
=$ $0$. We do not require $(\epsilon _{t},x_{t})$ to be mutually
independent, but use the second order condition $P(E[\epsilon
_{t}^{2}|x_{t}] $ $=$ $\sigma ^{2})$ $=$ $1$ for some finite $\sigma ^{2}$ $%
> $ $0$ to focus ideas. Thus $(\delta _{0},\theta _{0})$ should be seen as
being possibly pseudo-true in the sense of \cite{Sawa1978} and \cite%
{White1982}. Heterogeneity, including (conditional) heteroscedasticity, is
straightforward to include, but omitted here due to space considerations.
See the appendix for all assumptions. We can test general hypotheses $%
H_{0}:\theta _{0}=\tilde{\theta}$ vs. $H_{1}:\theta _{0}\neq \tilde{\theta}$%
, but testing $\theta _{0}$ $=$ $0$ naturally saves notation.

Although $k_{\theta ,n}$ $>>$ $n$, we do not require a sparsity condition
common in the high dimensional estimation literature 
\citep[see,
e.g.,][]{Hastie_etal2009,ZhangZhang2014,Dezeure_etal_2017}. There are many
uses for (\ref{model}), including panels with many fixed effects, estimation
with many instrumental variables, and linear projections of nonlinear
response, e.g. series expansions from flexible functional forms in machine
learning %
\citep[e.g.][]{Andrews1991,RoystonAltman1994,BelloniChernFern2017,CattaneoJanssonNewey2018,GautierRose2021}%
.

There is a large literature on high dimensional linear regression with an
increasing number of covariates; consult \cite{ZhangZhang2014}, \cite%
{Dezeure_etal_2017}, \cite{CattaneoJanssonNewey2018} and \cite{LiMuller2021}%
, amongst others. This literature concerns estimation for high dimensional
models, with inference \textit{post}-estimation. We, however, approach
inference on $\theta _{0}$ \textit{pre}-estimation, in the sense that (\ref%
{model}) is not estimated, hence a penalized or regularized estimator like
(de-biased) Lasso, a projection method, or partialing out parameters are not
required. Instead, we operate on many low dimension regression models, each
with one key parameter $\theta _{0,i}$ to be tested, and any other nuisance
parameters $\delta _{0}$ to be estimated. Parsimony improves estimation of
each $\theta _{0,i}$ when $k_{\theta ,n}$ is large, and allows us to
sidestep sparsity considerations under $H_{1}$. Indeed, pre-estimation
itself allows one to test a sparsity assumption based on a target set of
parameters.

Our\ test approach is as follows. Write $\beta _{(i)}$ $\equiv $ $[\delta
^{\prime },\theta _{i}]^{\prime }$, and construct parsimonious models:%
\begin{equation}
y_{t}=\delta _{(i)}^{\ast \prime }x_{\delta ,t}+\theta _{i}^{\ast }x_{\theta
,i,t}+v_{(i),t}=\beta _{(i)}^{\ast \prime }x_{(i),t}+v_{(i),t},\text{ }%
i=1,...,k_{\theta ,n},  \label{linear_parsim}
\end{equation}%
where (pseudo-true) $\beta _{(i)}^{\ast }$ for each $i$ satisfy $E[(y_{t}$ $%
- $ $\beta _{(i)}^{\ast \prime }x_{(i),t})x_{(i),t}]$ $=$ $0$. Thus the $%
i^{th} $ model corresponds to a single unique regressor $x_{\theta ,i,t}$
from $x_{\theta ,t}$ $=$ $[x_{\theta ,i,t}]_{i=1}^{k_{\theta ,n}}$, and
nuisance regressor subset $x_{\delta ,t}$. The methods proposed here can be
directly extended to more general settings, including nonlinear regression
models with non-identically distributed observations, conditional moment
models, random choice and quantile regression to name a few.

A key result shows that under mild conditions, $\theta ^{\ast }$ $\equiv $ $%
[\theta _{i}^{\ast }]_{i=1}^{k_{\theta ,n}}$ $=$ $0$\ \textit{if and only if}
$\theta _{0}$ $=$ $0$, thus each nuisance parameter $\delta _{(i)}^{\ast }$
is identically the true $\delta _{0}$. Obviously we cannot generally
identify $\theta _{0}$ (or $\delta _{0}$) under the alternative, but the set
of models in (\ref{linear_parsim}) can identify whether $\theta _{0}$ $=$ $0$
is true or not. Thus (\ref{linear_parsim}) represents the set of \textit{%
least parameterized} models, and therefore the fewest models ($k_{\theta ,n}$%
), required to identify whether the null hypothesis is correct, leading to
an efficiency gain.

Now let $\hat{\beta}_{(i)}$ $\equiv $ $[\hat{\delta}_{(i)}^{\prime },\hat{%
\theta}_{i}]^{\prime }$ minimize the parsimonious least squares loss $%
1/n\sum_{t=1}^{n}(y_{t}$ $-$ $\beta _{(i)}^{\prime }x_{(i),t})^{2}$. The
test statistic is the normalized maximum weighted $\hat{\theta}_{i}$:%
\begin{equation*}
\mathcal{T}_{n}=\max_{1\leq i\leq k_{\theta ,n}}\left\vert \sqrt{n}\mathcal{W%
}_{n,i}\hat{\theta}_{i}\right\vert ,
\end{equation*}%
where $\{\mathcal{W}_{n,i}\}_{n\geq 1}$ are sequences of possibly stochastic
weights, $\mathcal{W}_{n,i}$ $>$ $0$ $a.s.$ for each $i$, with non-random
(probability) limits $\mathcal{W}_{i}$ $\in $ $(0,\infty )$. The max-test
rejects $H_{0}$ at level $\alpha $ when an appropriate p-value approximation 
$\hat{p}_{n}$\ satisfies $\hat{p}_{n}$ $<$ $\alpha $. The weights $\mathcal{W%
}_{n,i}$\ allow for a variety of test statistics, including $\max_{1\leq
i\leq k_{\theta ,n}}|\sqrt{n}\hat{\theta}_{i}|$, or a max-t-statistic.
Bootstrap inference for the max-statistic does not require a covariance
matrix inversion like Wald and LM tests, and operates like a shrinkage
estimator by using only the most relevant (weighted) estimator. Our high
dimensional parametric wild bootstrap theory exploits multiplier bootstrap
theory in \cite{Chernozhukov_etal2013}, and the concept of \textit{weak
convergence in probability} in \cite{GineZinn1990}. The max-test is
asymptotically correctly sized, and consistent because $k_{\theta ,n}$ $%
\rightarrow $ $\infty $. Bootstrapping in high dimension has a rich history,
e.g. \cite{BickelFreedman1983}, \cite{Mammen1993}, %
\citet{Chernozhukov_etal2013,Chernozhukov_etal2017} and \cite%
{Dezeure_etal_2017} and their references.

We detail four covariate cases in order to bound $k_{\theta ,n}$, yielding $%
\ln (k_{\theta ,n})$ $=$ $o(n^{1/7})$ under $(i)$ boundedness, $(ii)$
sub-Gaussianicity, or $(iii)$ sub-exponentiality; and $(iv)$ under $\mathcal{%
L}_{p}$-boundedness $o((n/\ln (n)^{2})^{p/8})$ where $p$ $\geq $ $4$.
Furthermore, the max-test has non-trivial power against a sequence of $\sqrt{%
n/\{\ln \left( k_{\theta ,n}\right) E[\max_{1\leq i\leq k_{\theta
,n}}|x_{(i),t}|^{4}]\}}$-local alternatives.

The dimension $k_{\delta }$ of $\delta _{0}$\ is assumed fixed and $%
k_{\delta }$ $<$ $n$ to focus ideas, but can in principle vary with $n$, or
be unbounded. In a more general framework, if $\delta _{0}$ is sparse then a
penalized estimator can be used to estimate $[\delta _{(i)}^{\ast },\theta
_{i}^{\ast }]$, like DBL, ridge or Dantzig. If $k_{\delta }$ $\rightarrow $ $%
\infty $ as $n$ $\rightarrow $ $\infty $ then the nuisance variables can be
partialed out, similar to \cite{CattaneoJanssonNewey2018}. If $\delta _{0}$
is a (infinite dimensional) function, e.g. $y_{t}$ $=$ $\delta
_{0}(x_{\delta ,t})$ $+$ $\theta _{0}^{\prime }x_{\theta ,t}$ $+$ $\epsilon
_{t}$, then $[\delta _{(i)}^{\ast },\theta _{i}^{\ast }]$ can be estimated
by sieves \citep[e.g.][]{ChenPouzo2015}. In our setting inference is on the
high dimensional $\theta _{0}$, contrary to Chernozhukov et al (%
\citeyear{Chernozhukov_etal2016_DML}) and \cite{CattaneoJanssonNewey2018}\
who test a fixed low dimension parameter after partialing out a high
dimensional term.

Conventional test statistics typically exhibit size distortions due both to
poor estimator sample properties, and from inverting a high dimension
covariance matrix estimator that may be a poor approximation of the small or
asymptotic variance. The challenge of high dimensional covariance matrix
estimation is well documented \citep[e.g.][]{ChenXuWu2015}. A bootstrap
method is therefore typically applied, but bootstrap tests may only have
size corrected power equal to the original test, which may be low under
parameter proliferation %
\citep[e.g.][]{DavidsonMacKinnon2006,GhyselsHillMotegi2020}. Indeed, the
bootstrapped Wald test when $k_{\theta ,n}$ $+$ $k_{\delta }$ $<$ $n$ is
undersized with low power in linear iid models when $k_{\delta }$ is even
mildly large (e.g. $k_{\delta }$ $=$ $10)$\footnote{%
This may carry over to (semi)nonparametric settings with fixed low dimension 
$\delta _{{}}$ and infinite dimensional function $\theta _{0}$. See \cite%
{ChenPouzo2015} for simulation evidence of sieve based bootstrapped t-tests.}%
, and potentially significantly under-sized with low power when $k_{\theta
,n}$ is large, becoming acute when there are nuisance parameters. Max-tests
yield better (typically sharp) size, while power improvements over the Wald
test can be sizable when deviations from the null are small.

Reduced dimension regression models in the high dimensional parametric
statistics and machine learning literatures are variously called \textit{%
marginal regression}, \textit{correlation learning}, and \textit{sure
screening} \citep[e.g.][]{Fan_Lv_2008,Genovese_etal_2012}, and includes 
\textit{canonical correlation analysis} [CCA] %
\citep[e.g.][]{McKeagueZhang2022}. \cite{McKeague_Qian_2015}, for example,
regress $y$ against one covariate $x_{i}$ at a time, $i$ $=$ $1,...,k$ $<$ $%
\infty $, and compute the most relevant index $\hat{k}$ $\equiv $ $\arg
\max_{1\leq i\leq k}|\hat{\theta}_{i}|$ where $\hat{\theta}_{i}$ $\equiv $ $%
\widehat{cov}(y,x_{i})/\widehat{var}(x_{i})$, and test $\tilde{H}_{0}$ $:$ $%
\theta _{0}$ $=$ $[cov(y,x_{i})/var(x_{i})]_{i=1}^{k}$ $=$ $0$ using $\hat{%
\theta}_{\hat{k}}$. Thus they work outside a general regression setting,
consider only a fixed number $k$ of possible covariates, do not consider a
subset of regressor predictability (e.g. $x_{\theta ,t}$) with other
nuisance parameters ($\delta _{0}$), and use linear dependence for defining
predictability. Related projections to low dimension settings are considered
in \cite{Buhlmann2013}, $\ $\cite{ZhangZhang2014}, and \cite%
{GhyselsHillMotegi2020}, amongst others. \cite{McKeagueZhang2022} operate on
the canonical correlation between bounded high dimensional variables $x_{t}$ 
$\in $ $\mathbb{R}^{p}$ and $y_{t}$ $\in $ $\mathbb{R}^{q}$. They propose
nonparametric max-type tests of the global null hypothesis that there are no
linear relationships between any subsets of $x_{t}$ and $y_{t}$. In
economics and finance settings, however, there will typically be nuisance
covariates to control for, and the potential for nonlinear relationships.

Compared to \cite{McKeague_Qian_2015} and \cite{McKeagueZhang2022}, we also
work pre-estimation, but within a high dimensional\ regression model that
allows for misspecification and nuisance parameters, and test the parameter
subset $H_{0}$ $:$ $\theta _{0}$ $=$ $0$ where $\theta _{0}$ may have a
diverging dimension. One could train $y_{t}$ on just the key subset $%
x_{\theta ,t}$, cf. \cite{McKeagueZhang2022}, but that neglects intermediary
relationships (or controls) with $x_{\delta ,t}$, and they impose
boundedness on both $x_{t}$ and $y_{t}$. We also do not impose sparsity
(under $H_{1}$), nor any distributional assumptions, other than bounding $%
k_{\theta ,n}$ under four covariate cases from boundedness to $\mathcal{L}%
_{p}$-boundedness with $p$ $\geq $ $4$. If there are no nuisance parameters (%
$k_{\delta }$ $=$ $0$) then our theory should allow for a high dimensional
expansion of McKeague and Qian's (\citeyear{McKeague_Qian_2015}) method, or
an extension allowing for partialled out nuisance parameter, but we leave
that arc for future consideration.

Moreover, our method has a straightforward extension to nonlinear models and
a broad class of extremum estimators so that predictability may not align
with linear dependence between $y_{t}$ and the $x_{i,t}^{\prime }s$ (as
opposed to correlations, cf. \cite{McKeague_Qian_2015} and \cite%
{McKeagueZhang2022}). Consider a nonlinear response $f(x_{t},\beta _{0})$
with known $f$ and parsimonious versions $f(x_{t},[\delta _{(i)}^{\ast
};0,...,\theta _{(i)}^{\ast },...0])$, where $[0,...,\theta
_{(i)},...0]^{\prime }$ is a $k_{\theta ,n}$ $\times $ $1$ zero vector with $%
\theta _{(i)}$\ in the $i^{th}$ row. Thus $\theta _{0}$ need not align
one-for-one with a regressor subset $x_{\cdot ,t}$. Examples abound in the
financial and macro-econometrics, neural network and machine learning
literatures, including switching models, flexible functional forms, and
basis expansions.

Forward stepwise regressions exploit parsimonious models with increasing
dimension. In the estimation literature low dimension, typically linear,
models are built by adding one covariate at a time %
\citep[e.g.][]{Amemiya1980,Stone1981,RoystonAltman1994,CastleQinReed2009}.
Stepwise regression, like regularized estimators, depend on a reliable
choice of tuning parameter (or stopping rule). Our approach is different
because we do not build up a model; rather, we use a different class of
parsimonious models that combines just one key covariate with the nuisance
covariates. Our models are true under the null, whereas stepwise methods
typically do not hinge on a presumed hypothesis. Classically, \cite{Kabe1963}
tests whether a low dimension subset $\theta _{0}$ $=$ $0$ with stepwise
regressions that partial out $\delta _{0}$, effectively using parsimonious
models. This idea has been treated many times since, e.g. \cite%
{CattaneoJanssonNewey2018} and Chernozhukov et al (%
\citeyear{Chernozhukov_etal2016_DML}). Our method does not partial out low
or high dimensional components, nor add covariates one at time based on
predictor power.

There are several differences between our method and the max-test in \cite%
{Dezeure_etal_2017}. First, they consider different bootstrap techniques
designed to allow for heteroscedastic errors. As discussed above, we focus
on homoscedastic errors given space constraints. Second, we do not estimate (%
\ref{model}) directly, hence we do not require shrinkage under either
hypothesis. Third, we allow for any covariate as long as it has a fourth
moment, and derive bounds on $k_{\theta ,n}$ under four covariate cases.
Fourth, we provide a class of weighted max-statistics covering $\max_{1\leq
i\leq k_{\theta ,n}}|\sqrt{n}\hat{\theta}_{i}|$ as well as a max-t-test.
Fifth, computation time for our bootstrapped p-value is \textit{significantly%
} faster than required for the bootstrapped node-wise DBL because we do not
require a regularized estimator.\footnote{%
Using Matlab (with coordinate descent and ADMM algorithms) via a SLURM
scheduler on the Longleaf cluster at UNC, with 128 workers on 1 node, and
setting $n$ $=$ $100$, $k_{\delta }$ $=$ $0$, $k_{\theta ,n}$ $=$ $200$, and 
$1000$ bootstrap samples, \textit{one} bootstrapped p-value under $H_{0}$
for the DBL max-statistic with $5$-fold cross-validation took $335$ seconds (%
$5.6$ minutes). Increase $k_{\theta ,n}$ to $480$, as in our simulation
study, yielded a computation time of $8109$ seconds ($135.2$ minutes). The
bootstrapped parsimonious max-test, by comparison, took only $4.17$ and $5.1$
seconds, respectively, yielding up to a $1600x$ computation time gain.
Increase $n$ to $250$ and set $k_{\theta ,n}$ $=$ $1144$, and DBL required $%
50.21$ hours, while our method took just $5.2$ seconds (roughly $35,000x$
faster). See Section \ref{sec:sim} for complete simulation details.\label%
{fn:matlab}} For further reading on post-estimation hypothesis testing for
regularized estimators under sparsity, see \cite{Liu_Luo_2014}, and see \cite%
{Cai_Guo_2017} and \cite{CaiGuoXia2023} for minimax optimality properties
and recent surveys.

The remainder of this paper proceeds as follows. Section \ref{sec:pvalue}
covers hypothesis identification and test statistic asymptotics, and Section %
\ref{sec:param_boot} presents bootstrap p-value theory. A Monte Carlo study
is presented in Section \ref{sec:sim}, and parting comments are left for
Section \ref{sec:conclude}. Assumptions, proofs of main results, and a
subset of simulation results are relegated to the appendix. The supplemental
material [SM] \cite{supp_mat_testmanyzeros} contains omitted proofs and all
simulation results.\medskip\ 

We assume all random variables exist on a complete measure space. $|x|$ $=$ $%
\sum_{i,j}|x_{i,j}|$ is the $l_{1}$-norm, $|x|_{2}$ $=$ $%
(\sum_{i,j}x_{i,j}^{2})^{1/2}$ is the Euclidean or $l_{2}$ norm, and $||A||$ 
$=$ $\max_{|\lambda |_{2}}\{|A\lambda |_{2}/|\lambda |_{2}\}$ is the
spectral norm for finite dimensional square matrices $A$ (and the Euclidean
norm for vectors). $||\cdot ||_{p}$ denotes the $L_{p}$-norm. $a.s.$ is 
\textit{almost surely}. $\boldsymbol{0}_{k}$ denotes a zero vector with
dimension $k$ $\geq $ $1$. Write $r$-vectors as $x$ $\equiv $ $%
[x_{i}]_{i=1}^{r}$. $[\cdot ]$ rounds to the nearest integer. $K$ $>$ $0$ is
non-random and finite, and may take different values in different places. $%
awp1$ = asymptotically with probability approaching one. We say $z$ has 
\textit{sub-exponential} distribution tails when $P(|z|$ $>$ $\varepsilon )$ 
$\leq $ $b\exp \{c\varepsilon \}$ for some $(b,c)$ $>$ $0$ and all $%
\varepsilon $ $>$ $0$ \citep[see,
e.g.,][Chap.2.7]{Vershynin2018}.

\section{Max-Test and p-Value Computation\label{sec:pvalue}}

We first show that use of the set of parsimonious models $y_{t}=\delta
_{(i)}^{\ast \prime }x_{\delta ,t}+\theta _{i}^{\ast }x_{\theta
,i,t}+v_{(i),t}$, $i$ $=$ $1,...,k_{\theta ,n}$, may be used to construct a
test of $\theta _{0}$ $=$ $\boldsymbol{0}_{k_{\theta ,n}}$. Collect all $%
\theta _{i}^{\ast }$ in (\ref{linear_parsim}) into $\theta ^{\ast }$ $\equiv 
$ $\left[ \theta _{i}^{\ast }\right] _{i=1}^{k_{\theta ,n}}$.

\begin{theorem}
\label{th:ident}$\theta _{0}$ $=$ $\boldsymbol{0}_{k_{\theta ,n}}$ if and
only if $\theta ^{\ast }$ $=$ $\boldsymbol{0}_{k_{\theta ,n}}$. Under $H_{0}$
$:$ $\theta _{0}$ $=$ $\boldsymbol{0}_{k_{\theta ,n}}$ hence $\delta
_{(i)}^{\ast }$ $=$ $\delta _{0}$ $\forall i$, and under $H_{1}$ $:$ $H_{1}$ 
$:$ $\theta _{0,i}$ $\neq $ $0$ for at least one $i$ $\in $ $%
\{1,...,k_{\theta ,n}\}$ there exists an $i$ such that $\theta _{i}^{\ast }$ 
$\neq $ $0$.
\end{theorem}

\subsection{Max-Test Gaussian Approximation\label{sec:max_test_dist}}

Now define parsimonious least squares loss gradient and Hessians,%
\begin{equation*}
\widehat{\mathcal{G}}_{(i)}\equiv -\frac{1}{n}%
\sum_{t=1}^{n}v_{(i),t}x_{(i),t}\text{, }\widehat{\mathcal{H}}_{(i)}\equiv 
\frac{1}{n}\sum_{t=1}^{n}x_{(i),t}x_{(i),t}^{\prime }\text{, \ }\mathcal{H}%
_{(i)}\equiv E\left[ x_{(i),t}x_{(i),t}^{\prime }\right] ,
\end{equation*}%
and define the usual least squares first order terms: $\mathcal{\hat{Z}}%
_{(i)}$ $\equiv $ $-\sqrt{n}\mathcal{H}_{(i)}^{-1}\widehat{\mathcal{G}}%
_{(i)} $ $=$ $\mathcal{H}_{(i)}^{-1}n^{-1/2}\sum_{t=1}^{n}v_{(i),t}x_{(i),t}$%
. Under $H_{0}$ of course $v_{(i),t}$ $=$ $\epsilon _{t}$ for each $i$\ by
Theorem \ref{th:ident}. By standard least squares arguments $0$ $=$ $%
\widehat{\mathcal{G}}_{(i)}$ $+$ $\widehat{\mathcal{H}}_{(i)}(\hat{\beta}%
_{(i)}$ $-$ $\beta _{(i)}^{\ast })$, hence%
\begin{equation*}
\sqrt{n}\left( \hat{\beta}_{(i)}-\beta _{(i)}^{\ast }\right) =-\sqrt{n}%
\mathcal{H}_{(i)}^{-1}\widehat{\mathcal{G}}_{(i)}-\left\{ \widehat{\mathcal{H%
}}_{(i)}^{-1}-\mathcal{H}_{(i)}^{-1}\right\} \sqrt{n}\widehat{\mathcal{G}}%
_{(i)}\equiv \mathcal{\hat{Z}}_{(i)}+\mathcal{\hat{R}}_{i}.
\end{equation*}

We want a high dimensional asymptotic Gaussian approximation $|\max_{1\leq
i\leq k_{\theta ,n}}|\sqrt{n}\mathcal{W}_{n,i}\hat{\theta}_{i}|$ $-$ $%
\max_{1\leq i\leq k_{\theta ,n}}|\mathcal{W}_{i}\boldsymbol{Z}_{(i)}||$ $%
\overset{p}{\rightarrow }$ $0$ where $\{\boldsymbol{Z}_{(i)}\}_{i\in \mathbb{%
N}}$ is Gaussian. Set all weights $\mathcal{W}_{n,i}$ $=$ $1$ to ease
notation for now. Then: 
\begin{eqnarray}
&&\left\vert \max_{1\leq i\leq k_{\theta ,n}}\left\vert \sqrt{n}\hat{\theta}%
_{i}\right\vert -\max_{1\leq i\leq k_{\theta ,n}}\left\vert \boldsymbol{Z}%
_{(i)}\right\vert \right\vert  \label{expand} \\
&&\text{ \ \ \ \ \ \ \ \ \ \ \ \ \ \ \ }\leq \left\vert \max_{1\leq i\leq
k_{\theta ,n}}\left\vert [\boldsymbol{0}_{k_{\delta }}^{\prime },1]\mathcal{%
\hat{Z}}_{(i)}\right\vert -\max_{1\leq i\leq k_{\theta ,n}}\left\vert 
\boldsymbol{Z}_{(i)}\right\vert \right\vert +\max_{1\leq i\leq k_{\theta
,n}}\left\vert [\boldsymbol{0}_{k_{\delta }}^{\prime },1]\mathcal{\hat{R}}%
_{i}\right\vert .  \notag
\end{eqnarray}%
We use a now standard Gaussian approximation theory under $H_{0}$ to bound $%
\max_{1\leq i\leq k_{\theta ,n}}|[\boldsymbol{0}_{k_{\delta }}^{\prime },1]%
\mathcal{\hat{Z}}_{(i)}|$ $-$ $\max_{1\leq i\leq k_{\theta ,n}}|\boldsymbol{Z%
}_{(i)}|$, and high dimensional concentration (like) inequalities to bound
\linebreak $\max_{1\leq i\leq k_{\theta ,n}}|[\boldsymbol{0}_{k_{\delta
}}^{\prime },1]\mathcal{\hat{R}}_{i}|$. If, however, we have stochastic
weights $\{\mathcal{W}_{n,i}\}$ then we must contend with $\max_{1\leq i\leq
k_{\theta ,n}}|\mathcal{W}_{n,i}$ $-$ $\mathcal{W}_{i}|$ by assumption or by
case, as we below.

Write $\mathcal{M}_{n}$ $\equiv $ $E[\max_{1\leq i\leq k_{\theta
,n}}|x_{(i),t}|^{4}]$. We have the following general result proving $%
\max_{1\leq i\leq k_{\theta ,n}}|\mathcal{\hat{R}}_{i}|$ $\overset{p}{%
\rightarrow }$ $0$ provided $\ln (k_{\theta ,n})$ $=$ $o(\sqrt{n}/\mathcal{M}%
_{n})$. Here and in the sequel we explore four broad covariate cases in
order to bound $\mathcal{M}_{n}$ (and thereby bound $k_{\theta ,n}$), in
each case uniformly in $i$: 
\begin{equation}
\begin{tabular}{llll}
$\boldsymbol{(i)}$ & bounded $x_{(i),t}$ & $\boldsymbol{(ii)}$ & 
sub-Gaussian $\left\vert x_{(i),t}\right\vert ^{4}$ \\ 
$\boldsymbol{(iii)}$ & sub-exponential $\left\vert x_{(i),t}\right\vert ^{4}$
& $\boldsymbol{(iv)}$ & $\mathcal{L}_{p}$-bounded $x_{(i),t}$, $p$ $\geq $ $%
4 $%
\end{tabular}
\label{case}
\end{equation}%
\noindent Examples of ($ii$) and ($iii$) are, respectively, $x_{(i),t}$ $=$ $%
\ln |y_{(i),t}|$ or $x_{(i),t}$ $=$ $|y_{(i),t}|^{1/4}$\ for sub-Gaussian $%
y_{(i),t}$; and $x_{(i),t}$ $=$ $|y_{(i),t}|^{1/2}$ for sub-Gaussian $%
y_{(i),t}$.

Following standard high dimensional concentration theory, cf. Nemirovski's
inequality, the thinner the tails then the greater the allowed upper bound
on $k_{\theta ,n}$. Assumption \ref{assum:smooth} is discussed in Appendix %
\ref{app:assum}. Proofs of lemmas are presented in SM.

\begin{lemma}
\label{lm:expansion}Let $H_{0}$ and Assumption \ref{assum:smooth} hold. Let
the weight sequences $\{\mathcal{W}_{n,i}\}_{n\in \mathbb{N}}$ satisfy $%
\max_{1\leq i\leq k_{\theta ,n}}|\mathcal{W}_{n,i}$ $-$ $\mathcal{W}_{i}$ $|$
$=O_{p}(\sqrt{\ln \left( k_{\theta ,n}\right) \mathcal{M}_{n}/n})$ for
non-stochastic $\mathcal{W}_{i}$ $\in $ $(0,\infty )$. We have:%
\begin{equation}
\max_{1\leq i\leq k_{\theta ,n}}\left\vert [\boldsymbol{0}_{k_{\delta
}}^{\prime },1]\mathcal{\hat{R}}_{i}\right\vert =O_{p}\left( \ln \left(
k_{\theta ,n}\right) \mathcal{M}_{n}/\sqrt{n}\right) .  \label{theta_approx}
\end{equation}%
Then $\max_{1\leq i\leq k_{\theta ,n}}|[\boldsymbol{0}_{k_{\delta }}^{\prime
},1]\mathcal{\hat{R}}_{i}|$ $\overset{p}{\rightarrow }$ $0$ if under
covariate case $(i)$ $\ln (k_{\theta ,n})$ $=$ $o(\sqrt{n})$; $(ii)$ $\ln
(k_{\theta ,n})$ $=$ $o(n^{1/3})$; $(iii)$ $\ln (k_{\theta ,n})$ $=$ $%
o(n^{1/4})$; or $(iv)$ $k_{\theta ,n}$ $=$ $o((n/\ln (n)^{2})^{p/8})$ where $%
p$ $\geq $ $4$.
\end{lemma}

\begin{remark}
\normalfont The weight limit $\mathcal{W}_{i}$\ is assumed non-random
because the Gaussian approximation theory below requires $\mathcal{W}_{i}%
\boldsymbol{Z}_{(i)}(\cdot )$ to be Gaussian.
\end{remark}

\begin{remark}
\label{rm:H}\normalfont The presence of the $l_{\infty }$ fourth moment $%
\mathcal{M}_{n}$ arises from dealing with $\max_{1\leq i\leq k_{\theta ,n}}$ 
$||\widehat{\mathcal{H}}_{(i)}^{-1}-\mathcal{H}_{(i)}^{-1}||$ in a high
dimensional least squares environment, cf. Nemirovski's moment inequality %
\citep[e.g.][Lemma 14.24]{BuhlmannVanDeGeer2011}.
\end{remark}

\begin{remark}
\normalfont The $k_{\theta ,n}$ bounds are more lenient than for the
Gaussian approximation, below, except under $\mathcal{L}_{p}$-boundedness.
Ultimately this is due to the relative lack of complexity that is required
to prove $\max_{1\leq i\leq k_{\theta ,n}}|[\boldsymbol{0}_{k_{\delta
}}^{\prime },1]\mathcal{\hat{R}}_{i}|$ $\overset{p}{\rightarrow }$ $0$ via
Nemirovski's inequality, compared to a Gaussian approximation theory for $%
\max_{1\leq i\leq k_{\theta ,n}}|[\boldsymbol{0}_{k_{\delta }}^{\prime },1]%
\mathcal{\hat{Z}}_{(i)}|$ \citep[cf.]{Chernozhukov_etal2013}. Thus, in
practice the imperative bounds on $k_{\theta ,n}$\ under cases ($i$)-($iii$)
appear in Theorem \ref{thm:max_dist} below.
\end{remark}

The weight approximation $\max_{1\leq i\leq k_{\theta ,n}}|\mathcal{W}_{n,i}$
$-$ $\mathcal{W}_{i}$ $|$ $=O_{p}(\sqrt{\ln \left( k_{\theta ,n}\right) 
\mathcal{M}_{n}/n})$ seems unavoidable, but is justified by noting that it
holds (trivially) with flat $\mathcal{W}_{n,i}$ $=$ $1$, or inverted
standard errors. The latter is formally stated as follows. Define residuals $%
\hat{v}_{(i),t}$ $=$ $y_{t}$ $-$ $\hat{\beta}_{(i)}^{\prime }x_{(i),t}$.
Squared standard errors for $\sqrt{n}\hat{\theta}_{i}$ are $\widehat{%
\mathcal{S}}_{(i)}^{2}$ $\equiv $ $[\widehat{\mathcal{H}}_{(i)}^{-1}]_{1,1}%
\mathcal{\hat{V}}_{(i),n}^{2}$ where $\mathcal{\hat{V}}_{(i),n}^{2}$ $=$ $%
1/n\sum_{t=1}^{n}\hat{v}_{(i),t}^{2}$, and the asymptotic value is $\mathcal{%
S}_{(i)}^{2}$ $\equiv $ $[\mathcal{H}_{(i)}^{-1}]_{1,1}E[v_{(i),t}^{2}]$.
Now write $\mathcal{W}_{n,i}$ $=$ $1/\widehat{\mathcal{S}}_{(i)}$ and $%
\mathcal{W}_{i}$ $=$ $1/\mathcal{S}_{(i)}$. The following result with the
usual Slutsky theorem argument yields in general $\max_{1\leq i\leq
k_{\theta ,n}}|\mathcal{W}_{n,i}$ $-$ $\mathcal{W}_{i}$ $|$ $=O_{p}(\sqrt{%
\ln \left( k_{\theta ,n}\right) \mathcal{M}_{n}/n})$.

\begin{lemma}
\label{lm:se}Let Assumption \ref{assum:smooth} hold. Then $\max_{1\leq i\leq
k_{\theta ,n}}|\widehat{\mathcal{S}}_{(i)}^{2}$ $-$ $\mathcal{S}_{(i)}^{2}|$ 
$=$ $O_{p}(\sqrt{\ln \left( k_{\theta ,n}\right) \mathcal{M}_{n}/n})$.
\end{lemma}

Now with $\mathcal{\hat{Z}}_{(i)}$ $=$ $-\sqrt{n}\mathcal{H}_{(i)}^{-1}%
\widehat{\mathcal{G}}_{(i)}$ we derive by case sequences $\{k_{\theta ,n}\}$%
, and scalar normal random variables $\boldsymbol{Z}_{(i)}(\lambda )$ $\sim $
$N(0,E[\epsilon _{t}^{2}]\lambda ^{\prime }\mathcal{H}_{(i)}^{-1}\lambda )$,
such that for each $\lambda $ $\in $ $\mathbb{R}^{k_{\delta }+1}$, $\lambda
^{\prime }\lambda $ $=$ $1$, the Kolmogorov distance 
\begin{equation}
\rho _{n}(\lambda )\equiv \sup_{z\geq 0}\left\vert P\left( \max_{1\leq i\leq
k_{\theta ,n}}\left\vert \lambda ^{\prime }\mathcal{\hat{Z}}%
_{(i)}\right\vert \leq z\right) -P\left( \max_{1\leq i\leq k_{\theta
,n}}\left\vert \boldsymbol{Z}_{(i)}(\lambda )\right\vert \leq z\right)
\right\vert \rightarrow 0.  \label{Gauss_apprx}
\end{equation}%
In terms of the linear contrasts $\lambda ^{\prime }\mathcal{\hat{Z}}_{(i)}$
we specifically care about $\lambda $ $=$ $[\boldsymbol{0}_{k_{\delta
}}^{\prime },1]^{\prime }$ in view of expansion (\ref{theta_approx}), hence
we only require a pointwise limit theory vis-\`{a}-vis $\lambda $. The
approximation does not require standardized $\lambda ^{\prime }\mathcal{\hat{%
Z}}_{(i)}$ and $\boldsymbol{Z}_{(i)}(\lambda )$ because we assume $\lambda
^{\prime }\mathcal{H}_{(i)}^{-1}\lambda $ lies in a compact subset of $%
(0,\infty )$ asymptotically uniformly in $i$ and $\lambda $. We work in the
now seminal setting of \cite{Chernozhukov_etal2013}, cf. \cite%
{Chernozhukov_etal2017}.

\begin{theorem}
\label{thm:max_dist}Let $H_{0}$ and Assumption \ref{assum:smooth} hold. Then 
$\rho _{n}(\lambda )$ $=$ $O(1/g(n))$ $\rightarrow $ $0$ for any slowly
varying $g(n)$ $\rightarrow $ $\infty $, where by the (\ref{case}) covariate
cases $(i)$-$(iii)$ $\ln (k_{\theta ,n})=o(n^{1/7})$, and $(iv)$ $k_{\theta
,n}$ $=$ $o(n^{p/2}/[g(n)\ln ^{2}\left( n\right) ]^{p})$ where $p$ $\geq $ $%
4 $.
\end{theorem}

\begin{remark}
\normalfont In their post-estimation, DBL max-test, \cite{Dezeure_etal_2017}
consider only covariate boundedness ($i$). They implicitly require slowly
varying $g(n)$ $\rightarrow $ $\infty $ in their similar Gaussian
approximation theory. Boundedness, however, need not hold in practice, and
therefore does not fully describe the scope of required bounds on $k_{\theta
,n}$. Although we also achieve $\ln (k_{\theta ,n})=o(n^{1/7})$\ under
sub-Gaussian and sub-exponential cases, under $\mathcal{L}_{p}$-boundedness
the range reduces to $k_{\theta ,n}$ $=$ $o(n^{p/2}/[g(n)\ln ^{2}\left(
n\right) ]^{p})$ where slowly varying $g(n)$ guides convergence of the
Kolmogorov distance $\rho _{n}(\lambda )$. By construction, therefore, $%
k_{\theta ,n}$ $=$ $o(n^{p/2-\iota })$ for infinitessimal $\iota $ $>$ $0$.
\end{remark}

\begin{remark}
\normalfont Imposing slow variation $g(n)$ $\rightarrow $ $\infty $
optimizes the upper bound on $k_{\theta ,n}$ under ($iv$). If a faster rate
of convergence $\rho _{n}(\lambda )$ $\rightarrow $ $0$ is desired, this
will both reduce, and complicate solving for, an upper bound on $k_{\theta
,n}$.
\end{remark}

\begin{remark}
\normalfont$k_{\theta ,n}$ cannot have an exponential upper bound if higher
moments do not exist. In this case, e.g., at least $k_{\theta ,n}$ $=$ $%
o\left( n^{2}/[g(n)\ln ^{2}\left( n\right) ]^{4}\right) $ under $\mathcal{L}%
_{4}$-boundedness. If all moments exist, however, but a moment generating
function near zero does not exist (e.g. log-normal: else we revert to cases (%
$ii$) or ($iii$)), then we only need $k_{\theta ,n}$ $=$ $o\left( n^{\varphi
}\right) $ for arbitrarily large $\varphi $ $>$ $0$.
\end{remark}

We now have the main result of this section. The statistic $\mathcal{T}_{n}$ 
$=$ $\max_{1\leq i\leq k_{\theta ,n}}|\sqrt{n}\mathcal{W}_{n,i}\hat{\theta}%
_{i}|$ can be well approximated by the max-Gaussian process $\{\max_{1\leq
i\leq k_{\theta ,n}}|\mathcal{W}_{i}\boldsymbol{Z}_{(i)}([\boldsymbol{0}%
_{k_{\delta }}^{\prime },1])|\}_{n\in \mathbb{N}}$. The bounds on $k_{\theta
,n}$\ from the Gaussian approximation Theorem \ref{thm:max_dist} dominate
the asymptotic approximation Lemma \ref{lm:expansion} bounds except under $%
\mathcal{L}_{p}$-boundedness for finite $p$.

\begin{theorem}
\label{th:max_theta_hat}Let Assumption \ref{assum:smooth} hold, and assume
for non-stochastic $\mathcal{W}_{i}$ $\in $ $(0,\infty )$, $\max_{i\in 
\mathbb{N}}\mathcal{W}_{i}$ $<$ $\infty $, that $\max_{1\leq i\leq k_{\theta
,n}}|\mathcal{W}_{n,i}$ $-$ $\mathcal{W}_{i}$ $|$ $=O_{p}(\sqrt{\ln \left(
k_{\theta ,n}\right) \mathcal{M}_{n}/n})$. Assume by the (\ref{case})
covariate cases $(i)$-$(iii)$ $\ln (k_{\theta ,n})$ $=$ $o(n^{1/7})$, and $%
(iv)$ $k_{\theta ,n}$ $=$ $o((n/\ln (n)^{2})^{p/8})$ where $p$ $\geq $ $4$.$%
\medskip $\newline
$a.$ Under $H_{0}$, for the Gaussian random variables $\boldsymbol{Z}%
_{(i)}(\lambda )$ $\sim $ $N(0,E[\epsilon _{t}^{2}]\lambda ^{\prime }%
\mathcal{H}_{(i)}^{-1}\lambda )$ in (\ref{Gauss_apprx}),%
\begin{equation}
\left\vert P\left( \max_{1\leq i\leq k_{\theta ,n}}\left\vert \sqrt{n}%
\mathcal{W}_{n,i}\hat{\theta}_{i}\right\vert \leq z\right) -P\left(
\max_{1\leq i\leq k_{\theta ,n}}\left\vert \mathcal{W}_{i}\boldsymbol{Z}%
_{(i)}([\boldsymbol{0}_{k_{\delta }}^{\prime },1])\right\vert \leq z\right)
\right\vert \rightarrow 0\text{ }\forall z\geq 0.  \label{PP0}
\end{equation}%
Furthermore, $\max_{1\leq i\leq k_{\theta ,n}}|\sqrt{n}\mathcal{W}_{n,i}\hat{%
\theta}_{i}|$ $\overset{d}{\rightarrow }$ $\max_{i\in \mathbb{N}}|\mathcal{W}%
_{i}\boldsymbol{Z}_{(i)}([\boldsymbol{0}_{k_{\delta }}^{\prime },1])|$%
.\medskip \newline
$b.$ Under $H_{1}$, $\max_{1\leq i\leq k_{\theta ,n}}|\sqrt{n}\mathcal{W}%
_{n,i}\hat{\theta}_{i}|$ $\rightarrow $ $\infty $ for any monotonic sequence
of positive integers $\{k_{\theta ,n}\}$, $k_{\theta ,n}$ $\rightarrow $ $%
\infty $.
\end{theorem}

\begin{remark}
\label{rm:EVT}\normalfont We can attempt to characterize the max-statistic $%
\mathcal{T}_{n}$ null limit sequence \linebreak $\{\max_{1\leq i\leq
k}\left\vert \mathcal{W}_{i}\boldsymbol{Z}_{(i)}([\boldsymbol{0}_{k_{\delta
}}^{\prime },1])\right\vert \}_{k\in \mathbb{N}}$ after re-scaling. In
Appendix D of the SM we present a simplified environment such that under $%
H_{0}$ for some positive non-random $b_{n}$ $\sim $ $\sqrt{2\ln (k_{\theta
,n})}$, we have $\sqrt{2\ln (k_{\theta ,n})/E[\epsilon _{t}^{2}]}(\mathcal{T}%
_{n}$ $-$ $b_{n})$ $\overset{d}{\rightarrow }$ $\mathfrak{G}$, a Gumbel law.
The next section, however, delivers a wild bootstrapped p-value
approximation that allows for nuisance parameters and arbitrary covariate
dependence, and naturally affords better small sample inference, \emph{%
without} requiring knowledge of the limit law of a sequence of standardized $%
\mathcal{T}_{n}$.
\end{remark}

\section{Parametric Bootstrap\label{sec:param_boot}}

We utilize a variant of the wild (multiplier) bootstrap for inference %
\citep[see][]{Bose1988,Liu1988,Goncalves_Kilian_2004}. Define the restricted
estimator $\hat{\beta}^{(0)}$ $\equiv $ $[\hat{\delta}^{(0)\prime },%
\boldsymbol{0}_{k_{\theta ,n}}^{\prime }]^{\prime }$ under $H_{0}$ $:$ $%
\theta _{0}$ $=$ $\boldsymbol{0}_{k_{\theta ,n}}$, where $\hat{\delta}^{(0)}$
minimizes the restricted least squares criterion $\sum_{t=1}^{n}(y_{t}$ $-$ $%
\delta ^{\prime }x_{\delta ,t})^{2}$. Define residuals $\epsilon
_{n,t}^{(0)} $ $\equiv $ $y_{t}$ $-$ $\hat{\delta}^{(0)\prime }x_{\delta ,t},
$ draw iid $\{\eta _{t}\}_{t=1}^{n}$ from $N(0,1)$, and generate $%
y_{n,t}^{\ast }$ $\equiv $ $\hat{\delta}^{(0)\prime }x_{\delta ,t}$ $+$ $%
\epsilon _{n,t}^{(0)}\eta _{t}$. Now construct regression models $%
y_{n,t}^{\ast }$ $=$ $\beta _{(i)}^{\prime }x_{(i),t}$ $+$ $v_{n,(i),t}$ for
each $i$ $=$ $1,...,k_{\theta ,n}$, and let $\widehat{\tilde{\beta}}_{(i)}$ $%
=$ $[\widehat{\tilde{\delta}}_{(i)}^{\prime },\widehat{\tilde{\theta}}%
_{i}]^{\prime }$ be the least squares estimator of $\beta _{(i)}$. The
bootstrapped test statistic is 
\begin{equation}
\mathcal{\tilde{T}}_{n}\equiv \max_{1\leq i\leq k_{\theta ,n}}\left\vert 
\sqrt{n}\mathcal{W}_{n,i}\widehat{\tilde{\theta}}_{i}\right\vert .
\label{T_boot}
\end{equation}%
Repeat the above steps $\mathcal{R}$ times, each time drawing a new iid
sequence $\{\eta _{j,t}\}_{t=1}^{n}$, $j$ $=$ $1...\mathcal{R}$. This
results in a sequence of bootstrapped estimators $\{\widehat{\tilde{\theta}}%
_{j}\}_{j=1}^{\mathcal{R}}$ and test statistics $\{\mathcal{\tilde{T}}%
_{n,j}\}_{j=1}^{\mathcal{R}}$ that are iid conditional on the sample $%
\{x_{t},y_{t}\}_{t=1}^{n}$. The approximate p-value is 
\begin{equation}
\tilde{p}_{n,\mathcal{R}}\equiv \frac{1}{\mathcal{R}}\sum_{j=1}^{\mathcal{R}%
}I(\mathcal{\tilde{T}}_{n,j}>\mathcal{T}_{n}).  \label{p_boot}
\end{equation}

The above algorithm varies from \cite{Goncalves_Kilian_2004} since they
operate on the \textit{unrestricted} estimator in order to generate
residuals. We impose the null when we estimate $\beta _{0}$, thus reducing
dimensionality without affecting asymptotics under either hypothesis.
Further, in view of the Lemma \ref{lm:expansion} asymptotic expansion, the
parametric wild bootstrap is asymptotically equivalent to a wild
(multiplier) bootstrap applied to the first order term $\mathcal{\hat{Z}}%
_{(i)}$. Simply draw iid $\{\eta _{t}\}_{t=1}^{n}$ from $N(0,1)$, compute $%
\widehat{\mathcal{\tilde{Z}}}_{(i)}$ $\equiv $ $\sqrt{n}\widehat{\mathcal{H}}%
_{(i)}^{-1}\widehat{\mathcal{\tilde{G}}}_{(i)}$ where $\widehat{\mathcal{%
\tilde{G}}}_{(i)}$ $\equiv $ $-1/n\sum_{t=1}^{n}\eta _{t}(y_{t}$ $-$ $\hat{%
\beta}_{(i)}^{\prime }x_{(i),t})x_{(i),t}$, and then $\mathcal{\tilde{T}}%
_{n} $ $\equiv $ $\max_{1\leq i\leq k_{\theta ,n}}|\sqrt{n}\mathcal{W}_{n,i}%
\widehat{\mathcal{H}}_{(i)}^{-1}\widehat{\mathcal{\tilde{G}}}_{(i)}|$.

Now write $\beta ^{(0)}$ $=$ $[\delta ^{(0)\prime },\boldsymbol{0}%
_{k_{\theta ,n}}^{\prime }]^{\prime }$\ and $\beta _{(i)}^{(0)}$ $=$ $%
[\delta ^{(0)\prime },0]^{\prime }$, where $\delta ^{(0)}$ minimizes $%
E[(y_{t}$ $-$ $\delta ^{\prime }x_{\delta ,t})^{2}]$ on compact $\mathcal{D}$
$\subset $ $\mathbb{R}^{k_{\delta }}$. Define the normalized bootstrapped
gradient under the null:%
\begin{equation}
\mathcal{\tilde{Z}}_{(i)}^{(0)}\equiv \left( E\left[ x_{\delta ,t}x_{\delta
,t}^{\prime }\right] \right) ^{-1}\frac{1}{\sqrt{n}}\sum_{t=1}^{n}\eta
_{t}\left( y_{t}-\delta ^{(0)\prime }x_{\delta ,t}\right) x_{\delta ,t}.
\label{Z_tilda_0}
\end{equation}%
We have the following bootstrap expansion. Recall $\mathcal{M}_{n}$ $\equiv $
$E[\max_{1\leq i\leq k_{\theta ,n}}|x_{(i),t}|^{4}]$.

The following Lemmas \ref{lm:b_boot_expansion} and \ref{lm:max_dist_boot}
and Theorem \ref{thm:boot_theta} offer improved bounds on $k_{\theta ,n}$
due to the iid Gaussian multiplier term, and vagaries of the proofs of
expansion and Gaussian approximation Lemmas \ref{lm:b_boot_expansion} and %
\ref{lm:max_dist_boot}.

\begin{lemma}
\label{lm:b_boot_expansion}Let Assumption \ref{assum:smooth} hold, and let
the weights satisfy $\max_{1\leq i\leq k_{\theta ,n}}|\mathcal{W}_{n,i}$ $-$ 
$\mathcal{W}_{i}$ $|$ $=O_{p}(\sqrt{\ln \left( k_{\theta ,n}\right) \mathcal{%
M}_{n}/n}))$ for non-stochastic $\mathcal{W}_{i}$ $\in $ $(0,\infty )$. Then:%
\begin{equation*}
\left\vert \max_{1\leq i\leq k_{\theta ,n}}\left\vert \sqrt{n}\mathcal{W}%
_{n,i}\widehat{\tilde{\theta}}_{i}\right\vert -\max_{1\leq i\leq k_{\theta
,n}}\left\vert \mathcal{W}_{i}[\boldsymbol{0}_{k_{\delta }}^{\prime },1]%
\mathcal{\tilde{Z}}_{(i)}^{(0)}\right\vert \right\vert =O_{p}\left( \ln
\left( k_{\theta ,n}\right) \mathcal{M}_{n}^{3/4}/\sqrt{n}\right) .
\end{equation*}%
Hence the $O_{p}(\cdot )$ term is $o_{p}(1)$ if\ by the (\ref{case})
covariate cases: $(i)$ $\ln (k_{\theta ,n})$ $=$ $o(n^{1/2})$; $(ii)$ $\ln
(k_{\theta ,n})$ $=$ $o(n^{4/11})$; $(iii)$ $\ln (k_{\theta ,n})$ $=$ $%
o(n^{2/7})$; and$\ (iv)$ $k_{\theta ,n}$ $=$ $o((n/[\ln \left( n\right)
]^{2})^{p/6})$, $p$ $\geq $ $4$.
\end{lemma}

\begin{remark}
\normalfont The upper bound $O_{p}(\ln \left( k_{\theta ,n}\right) \mathcal{M%
}_{n}^{3/4}/\sqrt{n})$ is smaller than the corresponding bound for $\hat{%
\theta}_{i}$ by a factor of $\mathcal{M}_{n}^{1/4}$ (hence convergence to
zero is faster), ultimately due to the Gaussian multiplier. This in turn
leads to greater upper bounds on $k_{\theta ,n}$ by case compared to Lemma %
\ref{lm:expansion}.
\end{remark}

In order to characterize a conditional high dimensional Gaussian
approximation for $\mathcal{\tilde{Z}}_{(i)}^{(0)}$, define $\mathcal{\bar{H}%
}_{(i)}^{(0)}$ $\equiv $ $E[x_{\delta ,t}x_{\delta ,t}^{\prime }]$ and%
\begin{equation}
\tilde{\sigma}_{(i)}^{2}(\lambda )\equiv \lambda ^{\prime }\left( \mathcal{%
\bar{H}}_{(i)}^{(0)}\right) ^{-1}E\left[ \left( y_{t}-\delta ^{(0)\prime
}x_{\delta ,t}\right) ^{2}x_{\delta ,t}x_{\delta ,t}^{\prime }\right] \left( 
\mathcal{\bar{H}}_{(i)}^{(0)}\right) ^{-1}\lambda .  \label{Vni0}
\end{equation}%
Notice $\tilde{\sigma}_{(i)}^{2}(\lambda )$ $=$ $E[\epsilon _{t}^{2}]\lambda
^{\prime }\mathcal{H}_{(i)}^{-1}\lambda $ under $H_{0}$ follows from Theorem %
\ref{th:ident}. As usual, we only care about $\lambda $ $=$ $[\boldsymbol{0}%
_{k_{\delta }}^{\prime },1]^{\prime }$. Let $\Rightarrow ^{p}$ denote 
\textit{weak convergence in probability} (Gin\'{e} and Zinn, %
\citeyear{GineZinn1990}: Section 3), a useful notion for multiplier
bootstrap asymptotics. We have the following conditional multiplier
bootstrap central limit theorem.

\begin{lemma}
\label{lm:max_dist_boot}Let Assumption \ref{assum:smooth} hold. Let $\{%
\boldsymbol{\tilde{Z}}_{(i)}(\lambda )\}_{i\in \mathbb{N}}$, $\boldsymbol{%
\tilde{Z}}_{(i)}(\lambda )$ $\sim $ $N(0,\tilde{\sigma}_{(i)}^{2}(\lambda ))$%
, be an independent copy of the Theorem \ref{th:max_theta_hat}.a null
distribution process $\{\boldsymbol{Z}_{(i)}(\lambda )\}_{i\in \mathbb{N}}$,
that is independent of the asymptotic draw $\{x_{t},y_{t}\}_{t=1}^{\infty }$%
. Let $\{k_{\theta ,n}\}$ satisfy by the (\ref{case}) covariate cases $(i)$-$%
(iii)$ $\ln \left( k_{\theta ,n}\right) $ $=$ $o(n^{1/3})$, and $(iv)$ $%
k_{\theta ,n}$ $=$ $o((n/\ln (n))^{p/4})$. Then $\sup_{z\geq
0}|P(\max_{1\leq i\leq k_{\theta ,n}}|\lambda ^{\prime }\mathcal{\tilde{Z}}%
_{(i)}^{(0)}|$ $\leq $ $z|\mathfrak{S}_{n})$ $-$ $P(\max_{1\leq i\leq
k_{\theta ,n}}|\boldsymbol{\tilde{Z}}_{(i)}(\lambda )|$ $\leq $ $z)|$ $%
\overset{p}{\rightarrow }$ $0.$ Furthermore $\max_{1\leq i\leq k_{\theta
,n}}|[\boldsymbol{0}_{k_{\delta }}^{\prime },1]\mathcal{\tilde{Z}}%
_{(i)}^{(0)}|$ $\Rightarrow ^{p}$ $\max_{i\in \mathbb{N}}|\boldsymbol{\tilde{%
Z}}_{(i)}([\boldsymbol{0}_{k_{\delta }}^{\prime },1]\mathfrak{)}|$.
\end{lemma}

Bootstrap asymptotic expansion Lemma \ref{lm:b_boot_expansion}, and
bootstrap Gaussian approximation Lemma \ref{lm:max_dist_boot}, coupled with $%
|\mathcal{W}_{n,i}$ $-$ $\mathcal{W}_{i}|_{\infty }$ $=$ $O_{p}(\sqrt{\ln
\left( k_{\theta ,n}\right) \mathcal{M}_{n}/n}))$ prove the following claim.
The bounds on $k_{\theta ,n}$ are merely the least yielded from the two
supporting lemmas.

\begin{theorem}
\label{thm:boot_theta}Let Assumption \ref{assum:smooth} hold, and for
non-stochastic $\mathcal{W}_{i}$ $\in $ $(0,\infty )$, $\max_{i\in \mathbb{N}%
}\mathcal{W}_{i}$ $<$ $\infty $, assume $\max_{1\leq i\leq k_{\theta ,n}}|%
\mathcal{W}_{n,i}$ $-$ $\mathcal{W}_{i}$ $|$ $=O_{p}(\sqrt{\ln \left(
k_{\theta ,n}\right) \mathcal{M}_{n}/n})$. Let $\{k_{\theta ,n}\}$ satisfy
by the (\ref{case}) covariate cases $(i)$ $\ln \left( k_{\theta ,n}\right) $ 
$=$ $o(n^{1/3})$; $(ii)$ $\ln \left( k_{\theta ,n}\right) $ $=$ $o(n^{1/3})$%
; $(iii)$ $\ln \left( k_{\theta ,n}\right) $ $=$ $o(n^{2/7})$; and $(iv)$ $%
k_{\theta ,n}$ $=$ $o((n/[\ln \left( n\right) ]^{2})^{p/6})$, $p$ $\geq $ $4$%
. Let $\{\boldsymbol{\tilde{Z}}_{(i)}(\lambda )\}_{i\in \mathbb{N}}$, $%
\boldsymbol{\tilde{Z}}_{(i)}(\lambda )$ $\sim $ $N(0,\tilde{\sigma}%
_{(i)}^{2}(\lambda ))$, be an independent copy of the Theorem \ref%
{th:max_theta_hat}.a null distribution process $\{\boldsymbol{Z}%
_{(i)}(\lambda )\}_{i\in \mathbb{N}}$, that is independent of the asymptotic
draw $\{x_{t},y_{t}\}_{t=1}^{\infty }$. Then $\max_{1\leq i\leq k_{\theta
,n}}|\sqrt{n}\mathcal{W}_{n,i}\widehat{\tilde{\theta}}_{i}|$ $\Rightarrow
^{p}$ $\max_{i\in \mathbb{N}}|\mathcal{W}_{i}\boldsymbol{\tilde{Z}}_{(i)}([%
\boldsymbol{0}_{k_{\delta }}^{\prime },1]\mathfrak{)}|$.
\end{theorem}

We now have the main result of this section: the approximate p-value $\tilde{%
p}_{n,\mathcal{R}_{n}}$ promotes a correctly sized and consistent test
asymptotically. Notice that while Lemmas \ref{lm:b_boot_expansion} and \ref%
{lm:max_dist_boot} offer improvements on the upper bound for $k_{\theta ,n}$%
, here we have the reduced bounds set by Theorem \ref{th:max_theta_hat}
because $\tilde{p}_{n,\mathcal{R}_{n}}$ is computed from the bootstrapped $%
\widehat{\tilde{\theta}}_{i}$ \textit{and} original $\hat{\theta}_{i}$.

\begin{theorem}
\label{th:p_value_boot}Let Assumption \ref{assum:smooth} hold, and for
non-stochastic $\mathcal{W}_{i}$ $\in $ $(0,\infty )$, $\max_{i\in \mathbb{N}%
}\mathcal{W}_{i}$ $<$ $\infty $, assume $\max_{1\leq i\leq k_{\theta ,n}}|%
\mathcal{W}_{n,i}$ $-$ $\mathcal{W}_{i}$ $|$ $=O_{p}(\sqrt{\ln \left(
k_{\theta ,n}\right) \mathcal{M}_{n}/n})$. Let $\{k_{\theta ,n}\}$ satisfy
by the (\ref{case}) covariate cases $(i)$-$(iii)$ $\ln (k_{\theta ,n})$ $=$ $%
o(n^{1/7})$, and $(iv)$ $k_{\theta ,n}$ $=$ $o((n/\ln (n)^{2})^{p/8})$ where 
$p$ $\geq $ $4$.\medskip \newline
$a.$ Under $H_{0}$, $P(\tilde{p}_{n,\mathcal{R}_{n}}$ $<$ $\alpha )$ $%
\rightarrow $ $\alpha $; and $b.$ Under $H_{1}$, $P(\tilde{p}_{n,\mathcal{R}%
_{n}}$ $<$ $\alpha )$ $\rightarrow $ $1$.
\end{theorem}

\begin{remark}
\label{rm:null_increase}\normalfont The test is consistent against a global
alternative despite the fact that $\mathcal{T}_{n}$ increases in magnitude
in $n$, even under $H_{0}$. The proof reveals $\tilde{p}_{n,\mathcal{R}_{n}}$
is asymptotically in probability equivalent to 
\begin{equation*}
P\left( \max_{1\leq i\leq k_{\theta ,n}}\left\vert \mathcal{W}_{i}%
\boldsymbol{\tilde{Z}}_{(i)}([\boldsymbol{0}_{k_{\delta }}^{\prime
},1])\right\vert >\max_{1\leq i\leq k_{\theta ,n}}\left\vert \sqrt{n}%
\mathcal{W}_{n,i}\left( \hat{\theta}_{i}-\theta _{0,i}\right) +\sqrt{n}%
\mathcal{W}_{i}\theta _{0,i}+o_{p}(1)\right\vert \right) ,
\end{equation*}%
where $\{\boldsymbol{\tilde{Z}}_{(i)}([\boldsymbol{0}_{k_{\delta }}^{\prime
},1])\}$ is an independent copy of zero-mean Gaussian $\{\boldsymbol{Z}%
_{(i)}(\lambda )\}_{i\in \mathbb{N}}$. The max-statistic is therefore $O_{p}(%
\sqrt{\ln \left( k_{\theta ,n}\right) \mathcal{M}_{n}}$ $+$ $\sqrt{n})$
under $H_{1}$, cf. Lemma \ref{lm:beta_hat_rate} in the appendix. Conversely,
standard Gaussian concentration bounds yield $\max_{1\leq i\leq k_{\theta
,n}}|\mathcal{W}_{i}\boldsymbol{\tilde{Z}}_{(i)}([\boldsymbol{0}_{k_{\delta
}}^{\prime },1])|$ $=$ $O_{p}(\sqrt{\ln \left( k_{\theta ,n}\right) })$, cf.
(\ref{PWZ}) in the appendix. This yields consistency $\tilde{p}_{n,\mathcal{R%
}_{n}}$ $\overset{p}{\rightarrow }$ $0$.
\end{remark}

Finally, the bootstrapped test has non-trivial power against a sequence of
local alternatives $H_{1}^{L}:\theta _{0}=c\sqrt{\ln \left( k_{\theta
,n}\right) \mathcal{M}_{n}/n}$ where $c$ $=$ $\left[ c_{i}\right]
_{i=1}^{k_{\theta ,n}}$. The local-to-null rate follows from global
max-asymptotics Lemma \ref{lm:beta_hat_rate}.c. The rate is generally less
than $\sqrt{n}$ because $(i)$ we need to deal with $\max_{1\leq i\leq
k_{\theta ,n}}||\widehat{\mathcal{H}}_{(i)}^{-1}-\mathcal{H}_{(i)}^{-1}||$
and hence the fourth moment envelope $\mathcal{M}_{n}$, cf. Remark \ref{rm:H}%
, where generally $\mathcal{M}_{n}$ $\rightarrow $ $\infty $ is possible
(e.g. for unbounded covariates); and $(ii)$ the test statistic under the
null is, $awp1$, the maximum of an increasing sequence of Gaussian random
variables $\max_{1\leq i\leq k_{\theta ,n}}|\mathcal{W}_{i}\boldsymbol{%
\tilde{Z}}_{(i)}([\boldsymbol{0}_{k_{\delta }}^{\prime },1]\mathfrak{)}|$,
and $\max_{1\leq i\leq k_{\theta ,n}}|\mathcal{W}_{i}\boldsymbol{\tilde{Z}}%
_{(i)}([\boldsymbol{0}_{k_{\delta }}^{\prime },1]\mathfrak{)}|$ $=$ $O_{p}(%
\sqrt{\ln (k_{\theta ,n})})$ as discussed above

\begin{theorem}
\label{th:local_power}Under the conditions of Theorem\ \ref{th:p_value_boot}%
, $\lim_{n\rightarrow \infty }P(\tilde{p}_{n,\mathcal{R}_{n}}$ $<$ $\alpha )$
$\nearrow $ $1$ under $H_{1}^{L}$\ monotonically as the maximum local drift $%
\max_{i\in \mathbb{N}}|c_{i}|$ $\nearrow $ $\ \infty $.
\end{theorem}

\begin{remark}
\normalfont Asymptotic local power therefore does not depend on the degree,
if any, of sparsity. By comparison, for example, \cite{Zhong_Chen_Xu_2013}
consider a high dimensional mean model $y_{i}$ $=$ $\theta _{0}$ $+$ $%
\epsilon _{i}$ and test $H_{0}$ $:$ $\theta _{0}$ $=$ $\boldsymbol{0}_{p}$,
where $y_{i}$ $\in $ $\mathbb{R}^{p}$. Under their alternative $p^{1-\xi }$
elements $\theta _{0,i}$ $\neq $ $0$ for some $\xi $ $\in $ $(1/2,1)$, where 
$\theta _{0,i}$ $=$ $\sqrt{r\ln (p)/n}$ for some $r$ $>$ $0$, and by
assumption $\ln (p)$ $=$ $o(n^{1/3})$, hence $\theta _{0,i}$ $=$ $%
o(1/n^{2/3})$. See also \cite{Arias-Castro_Candes_Plan_2011} and their
references. Our max-test is consistent against local deviations $\theta _{0}$
$=$ $c(\ln (k_{\theta ,n})/n)^{\zeta }$ for any $\zeta $ $\in $ $[0,1/2)$
and $c$ $\neq $ $\boldsymbol{0}_{k_{\theta ,n}}$, and has non-trivial power
when $\zeta $ $=$ $1/2$ and $c$ $\neq $ $\boldsymbol{0}_{k_{\theta ,n}}$.
Further, we allow for any positive number of elements $\theta _{0,i}$ $\neq $
$0$ under $H_{1}$ or $H_{1}^{L}$ (i.e. $1$ $\leq $ $\sum_{i=0}^{k_{\theta
,n}}I(\theta _{0,i}$ $\neq $ $0)$ $\leq $ $k_{\theta ,n}$). Thus, in the
simple setting of \cite{Zhong_Chen_Xu_2013}, alternatives can be closer to
the null than allowed here in a max-test setting, but sparsity is enforced.
Here, any degree of (non)sparsity is permitted under $H_{1}$, (as long as
some $\theta _{0,i}$ $\neq $ $\ 0$) or $H_{1}^{L}$ (provided $c_{i}$ $\neq $ 
$\ 0$ for some $i$).
\end{remark}

\section{Monte Carlo Experiments\label{sec:sim}}

We use $1000$ independently drawn samples of sizes $n$ $\in $ $%
\{100,250,500\}$. We perform two max-tests: one with flat $\mathcal{W}_{n,i}$
$=$ $1$, and the other $\mathcal{W}_{n,i}$ is equal to the inverted least
squares standard error of $\hat{\theta}_{n,i}$, resulting in the \textit{%
max-test} and \textit{max-t-test}. Other tests are discussed below. Unless
otherwise noted, all bootstrapped p-values are based on $\mathcal{R}$ $=$ $%
1000$ independently drawn samples from $N(0,1)$. Significance levels are $%
\alpha $ $\in $ $\{.01,.05,.10\}$.

The following presents our benchmark design. Robustness checks are performed
in Appendix E of the supplemental material.

\subsection{Design and Tests of Zero Restrictions}

\subsubsection{Set-Up}

The DGP is $y_{t}$ $=$ $\delta _{0}^{\prime }x_{\delta ,t}$ $+$ $\theta
_{0}^{\prime }x_{\theta ,t}$ $+$ $\epsilon _{t}$ where $\epsilon _{t}$ is
iid standard normal. The total number of covariates is $k$ $\equiv $ $%
k_{\delta }$ $+$ $k_{\theta ,n}$. We consider two opposing nuisance
parameter cases $k_{\delta }$ $\in $ $\{0,10\}$, and for $x_{\theta ,t}$ we
set $k_{\theta ,n}$ $\in $ $\{35,k_{1}(n),k_{2}(n)\}$. We use $k_{1}(n)$ $%
\equiv $ $[\exp \{3.2n^{1/7-\iota }\}]$ with $\iota $ $=$ $10^{-10}$, which
is valid asymptotically when $x_{t}$ is bounded (cf. Theorem \ref%
{th:p_value_boot}). We also use $k_{2}(n)$ $=$ $[.02n^{2}]$, which is valid
asymptotically in any case here.\footnote{%
In our design $x_{t}$ is bounded or Gaussian and is therefore $\mathcal{L}%
_{p}$-bounded for any $p$ $>$ $0$.} Thus in this experiment $k_{\theta ,100}$
$\in $ $\{35,200,482\}$, $k_{\theta ,250}$ $\in $ $\{35,1144,1250\}$ and $%
k_{\theta ,500}$ $\in $ $\{35,2381,5000\}$.

We considered nine total covariate cases $[x_{\delta ,t},x_{\theta ,t}]$
(three bound cases and three dependency cases). In three broad scenarios $%
x_{i,t}$ are random draws from a truncated standard normal distribution,
denoted $\mathcal{TN}(0,1,U)$, such that $-U$ $\leq $ $x_{i,t}$ $\leq $ $U$ $%
a.s.$, with bounds $U_{1}$ $=$ $2.5$, $U_{2}$ $=$ $10^{10}$, or $U_{3}$ $=$ $%
\infty $ (i.e. $x_{i,t}$ $\sim $ $\mathcal{N}(0,1)$). Only upper bound $%
U_{1} $ has an impact since under $U_{2}$ the bound was never surpassed in
our experiment. The theory under boundedness, however, allows for any $%
k_{\theta ,n}$ provided $\ln (k_{\theta ,n})$ $=$ $o(n^{1/7})$.

Under these three scenarios we then consider three dependence types for $%
[x_{\delta ,t},x_{\theta ,t}]$: first, $[x_{\delta ,t},x_{\theta ,t}]$ are
serially and mutually independent standard normals; second, they are
block-wise dependent normals ($x_{\delta ,t}$ and $x_{\theta ,t}$ are
mutually independent); and third, they are within and across block dependent
normals.

Type three covariates are drawn as follows. Combine $x_{t}$ $\equiv $ $%
[x_{\delta ,t}^{\prime },x_{\theta ,t}^{\prime }]^{\prime }\in $ $\mathbb{R}%
^{k}$, and let $(w_{t},v_{t})$\ be mutually independent draws $\mathcal{TN}%
(0,I_{k})$. The regressors are $x_{t}$ $=$ $Aw_{t}$ $+$ $v_{t}$. We randomly
draw each element of $A\in $ $\mathbb{R}^{k}$ from a uniform distribution on 
$[-1,1]$. If $A$ does not have full column rank then we add a randomly drawn 
$\iota $ from $[0,1]$ to each diagonal component (in this study every draw $%
A $ had full column rank).

Not surprisingly all tests perform better when signals are less entangled.
Since our $U_{2}$-bounded and $U_{3}$-unbounded designs lead to identical
results, we therefore only report results for $U_{1}$-bounded and $U_{3}$%
-unbounded $x_{t}$ under dependence type three (within and across block
dependent).

We fix $\delta $ $=$ $\mathbf{1}_{k_{\delta }}$, a $k_{\delta }$ $\times $ $%
1 $ vector of ones. The benchmark models are as follows. Under the null $%
\theta _{0}$ $=$ $[\theta _{0,i}]_{i=1}^{k_{\theta ,n}}$ $=$ $\boldsymbol{0}%
_{k_{\theta ,n}}$. The alternatives are Alt($i$) $\theta _{0,1}$ $=$ $.0015$
with $[\theta _{0,i}]_{i=2}^{k_{\theta ,n}}$ $=$ $\boldsymbol{0}_{k_{\theta
,n}-1}$; Alt$(ii)$ $\theta _{0,i}$ $=$ $.0002i/k_{\theta ,n}$ for $i$ $=$ $%
1,...,k_{\theta ,n}$; and Alt$(iii)$ $\theta _{0,i}$ $=$ $.00015$ for $i$ $=$
$1,...,10$. See Table \ref{tbl_h1} for reference. Thus, under Alt$(i)$ only
one element deviates from, but is close to, zero. Alts $(ii)$ and $(iii)$
have smaller values, but for several or all $\theta _{0,i}$.

\begin{table}[h]
\caption{Alternative Models}
\label{tbl_h1}
\begin{center}
\begin{tabular}{ll}
Alt$(i):$ & $\theta _{0,1}$ $=$ $.0015$ and $[\theta
_{0,i}]_{i=2}^{k_{\theta ,n}}$ $=$ $0$ \\ 
Alt$(ii):$ & $\theta _{0,i}$ $=$ $.0002\times i/k_{\theta ,n}$ for $i$ $=$ $%
1,...,k_{\theta ,n}$ \\ 
Alt$(iii):$ & $\theta _{0,i}$ $=$ $.00015$ for $i$ $=$ $1,...,10$%
\end{tabular}%
\end{center}
\end{table}

\subsubsection{Parsimonious Max-Tests}

We estimate $k_{\theta ,n}$ parsimonious models $y_{t}$ $=$ $\delta
_{(i)}^{\ast \prime }x_{\delta ,t}$ $+$ $\theta _{i}^{\ast }x_{\theta ,i,t}$ 
$+$ $v_{(i),t}$ by least squares. Denote by $\mathcal{T}_{n}$ the resulting
max-test or max-t-test statistic.\medskip \textbf{\ }The bootstrapped test
statistic $\mathcal{\tilde{T}}_{n}$ and p-value $\tilde{p}_{n,\mathcal{R}}$
are computed as in (\ref{T_boot}) and (\ref{p_boot}). We reject $H_{0}$ when 
$\tilde{p}_{n,\mathcal{R}}$ $<$ $\alpha $.

\subsubsection{De-biased Lasso Max-Tests}

We follow \cite{Dezeure_etal_2017}, cf. \cite{ZhangZhang2014}, and estimate (%
\ref{model}) by DBL, using $5$-fold cross-validation in order to select the
tuning parameter. \cite{Dezeure_etal_2017} propose a (post-estimation)
max-t-test for a parameter subset, and various bootstrap techniques. We
perform both max- and max-t-tests for direct comparisons with the
parsimonious max-tests. The max-t-test uses the standard error in 
\citet[eq.
(4)]{Dezeure_etal_2017}, de-biasing is performed node-wise only for $\theta $%
, and p-values are approximated by parametric wild bootstrap as above.

Computation time for bootstrapped DBL is prohibitive. In Footnote \ref%
{fn:matlab} we saw our method is up to $1600x$ faster when $n$ $=$ $100$ and 
$35,000x$ faster when $n$ $=$ $250$, based on using $128$ workers in a
parallel processing environment. The speed difference, however, clouds just
how slow bootstrapped DBL is: just one p-value under $n$ $=$ $100$ with $%
k_{\theta ,n}$ $=$ $480$ took over $2$ hours, and over $50$ hours when $n$ $%
= $ $250$ and $k_{\theta ,n}$ $=$ $1144$. Scaling up to $1000$ samples is
obviously not currently attractive. We therefore ran a limited experiment
for DBL using only $n$ $=$ $100$ with $250$ samples, and $250$\ bootstrap
samples.

\subsubsection{Wald Tests}

If $k_{\delta }$ $+$ $k_{\theta ,n}$ $<$ $n$ then we also perform asymptotic
and bootstrapped Wald tests, and asymptotic and bootstrapped normalized Wald
tests. This will highlight in non-high dimensional settings the advantages
of a max-test over a Wald test. We estimate $y_{t}$ $=$ $\delta _{0}^{\prime
}x_{\delta ,t}$ $+$ $\theta _{0}^{\prime }x_{\theta ,t}$ $+$ $\epsilon _{t}$
by least squares. The asymptotic Wald test, based on the $\chi
^{2}(k_{\theta ,n})$ distribution, leads to large empirical size distortions
and is therefore not reported here. The bootstrapped Wald test is based on a
parametric wild bootstrap.

The normalized Wald statistic is $\mathcal{W}_{n}^{s}$ $\equiv $ $(\mathcal{W%
}_{n}$ $-$ $k_{\theta ,n})/\sqrt{2k_{\theta ,n}}$. Under the null $\mathcal{W%
}_{n}^{s}$ $\overset{d}{\rightarrow }$ $N(0,1)$ as $k_{\theta ,n}$ $%
\rightarrow $ $\infty $, and as long as $k_{\theta ,n}/n$ $\rightarrow $ $0$
then $\mathcal{W}_{n}^{s}$ $\overset{p}{\rightarrow }$ $\infty $ under $%
H_{1} $. This asymptotic one-sided test rejects the null when $\mathcal{W}%
_{n}^{s}$ $>$ $Z_{\alpha }$, where $Z_{\alpha }$ is the standard normal
upper tail $\alpha $-level critical value. This test yields highly distorted
empirical sizes and is therefore not reported. The bootstrapped test uses $%
\mathcal{W}_{n,i}^{s\ast }$ $\equiv $ $(\mathcal{W}_{n,i}^{\ast }$ $-$ $%
k_{\theta ,n})/\sqrt{2k_{\theta ,n}}$ with p-value approximation $p_{n}$ $=$ 
$1/\mathcal{R}\sum_{i=1}^{\mathcal{R}}I(\mathcal{W}_{n}^{s}$ $>$ $\mathcal{W}%
_{n,i}^{s\ast })$. Trivially $\mathcal{W}_{n}^{s}$ $>$ $\mathcal{W}%
_{n,i}^{s\ast }$ \textit{if and only if} $\mathcal{W}_{n,i}^{\ast }$ $>$ $%
\mathcal{W}_{n}$, hence we only discuss the bootstrapped Wald test.

\subsection{Simulation Results}

\subsubsection{Benchmark Results\label{sec:sim_bench}}

We report rejection frequencies in Tables \ref{table:max_wald_rejH0H11} and %
\ref{table:max_wald_rejH12H13} when nuisance terms $k_{\delta }$ $=$ $0$,
with unbounded $x_{t}$. See Appendix F of SM for complete simulation results
(bounded and unbounded $x_{t}$; $k_{\delta }$ $\geq $ $0$).

\paragraph{Empirical Size}

Let $p$-max and $dbl$-max denote parsimonious and DBL max-tests. Test
results under the three covariate cases are similar, although all tests
perform slightly less well under correlated regressors within and across
blocks (case three). We therefore only report and comment on the latter
results.

The $p$-max tests typically lead to qualitatively similar results with
empirical size close to nominal size (see Table \ref{table:max_wald_rejH0H11}%
). The tests perform roughly the same whether covariates are bounded or not.
In the presence of nuisance parameters $k_{\theta ,n}$ $=$ $10$ the tests
are slightly over-sized when $n$ $=$ $100$ with improvements as $n$
increases (the tests are comparable when $n$ $\geq $ $250$).

The $dbl$-max tests yield competitive empirical size when $k_{\delta }$ $=$ $%
0$, although the $dbl$-max test is generally over-sized compared to the $p$%
-max tests, especially when $k_{\theta ,n}$ is large. Furthermore, the $dbl$%
-max t-test tends to be highly under-sized when there are nuisance
parameters ($k_{\delta }$ $=$ $10$). The latter naturally implies
comparatively low power, discussed below. Overall, across cases the $p$%
-tests yield better empirical size than the $dbl$-max tests.

In the case $k_{\delta }$ $+$ $k_{\theta ,n}$ $<$ $n$, the asymptotic Wald
test and asymptotic normalized Wald test are generally severely over-sized
due to the magnitude of $k_{\delta }$ $+$ $k_{\theta ,n}$ (the latter test
is not shown). The bootstrapped Wald test is strongly under-sized when $n$
is small, but also when $k_{\delta }$ $=$ $10$ and/or when $k_{\theta ,n}$
is large.

\paragraph{Empirical Power}

See Tables \ref{table:max_wald_rejH0H11} and \ref{table:max_wald_rejH12H13}.
The bootstrapped Wald test in a low dimension setting yields lower empirical
power than the $p$-max tests in every case. The asymptotic Wald test yields
size corrected power that is nearly zero. We now focus on the max tests.

\subparagraph{Alt$(i)$}

Only one parameter exhibits a small deviation from the null: $\theta _{0,1}$ 
$=$ $.0015$. The difficulty in detecting such a deviation is apparent at $n$ 
$=$ $100$ when nuisance parameters are present and/or when $k_{\theta ,n}$
is small. The latter suggests a noise-signal effect: when $k_{\theta ,n}>>n$
the impact of the $k_{\delta }$ $=$ $10$ nuisance parameters is negligible,
but the impact can be large under any alternative studied here when $%
k_{\theta ,n}$ is small. For example, the $p$-max-tests yield (nearly) $%
100\% $ empirical power for sample sizes ($100$) $250$ and $500$ when $%
k_{\delta }$ $=$ $0$, or when $k_{\delta }$ $=$ $10$ and $k_{\theta ,n}>>n$.

If $n$ $=$ $100$, $k_{\theta ,n}$ $=$ $200$, $k_{\delta }$ $=$ $10$ and
covariates are bounded then power drops to $(.941,.975,.988)$ for $p$-max-t
at sizes $(.01,.05,.10)$, and when covariates are unbounded power is
comparable at $(.974,.993,.997)$. At $n=250$ with \ $k_{\theta ,n}$ $>>$ $n$
power rises to $100\%$ for bounded and unbounded covariates. Overall, $p$%
-max-test performance is not systematically sensitive to whether the
covariates are bounded.

The $dbl$-max-tests are generally dominated by the $p$-max tests, in many
cases significantly so. Consider $n$ $=$ $100$, $k_{\theta ,n}$ $=$ $200$, $%
k_{\delta }$ $=$ $10$ with bounded covariates: empirical power is only $%
(.002,.002,.003)$ for the $dbl$-max-t-test, compared to $(.941,.975,.988)$
for the $p$-max-t-test. If covariates are unbounded and there are no
nuisance parameters, then $dbl$-max-t-test power improves to $%
(.132,.188,.220)$, compared to $p$-max-t-test power $(1.00,1.00,1.00)$.

The $dbl$-max-test yields generally low power even when $k_{\delta }$ $=$ $0$%
: under $n$ $=$ $100$, $k_{\theta ,n}$ $=$ $482$ with unbounded covariates, $%
dbl$-max power is $(036,.180,.332)$ compared to $(.804,.900,.904)$ for $dbl$%
-max-t. However, when $k_{\theta ,n}$ $=$ $200$ power for both tests is
under $22\%$ for any nominal size compared to $100\%$ across sizes for $p$%
-max and $p$-max-t.

Considering $dbl$-max has sub-par power, in the remaining discussion we
focus only on $dbl$-max-t.

\subparagraph{Alt$(ii)$}

In this case all parameters deviate from zero monotonically between $%
(0,.0002]$. The $p$-max test yields $100\%$ empirical power when $k_{\delta
} $ $=$ $0$, or when $k_{\delta }$ $=$ $10$ and $k_{\theta ,n}$ $>>$ $n$,
irrespective of covariate (un)boundedness, again suggesting a signal/noise
effect.

The $dbl$-max tests, by comparison, are incapable of detecting such a signal
when nuisance parameters are present and $n$ $=$ $100$, yielding lower
rejection rates than nominal size. That power is trivial may be an artifact
of the design, and the requirements for DBL:\ $(i)$ We do not need to impose
sparsity because our method does not require it. Yet $(ii)$ DBL requires
sparsity \citep[e.g.][(B.2)]{Dezeure_etal_2017}: DBL is both picking up the
strong nuisance parameter values \textit{and} effectively setting $\hat{%
\theta}_{i}$ $\approx $ $0$. And $(iii)$ bias follows under sparsity
failure, in this case resulting in rejection rates below size.

The $dbl$-max tests yield decent empirical power when nuisance parameters
are not present, although consistently below, or far below, $p$-max tests.
This applies whether covariates are bounded 
\citep[as
in][]{Dezeure_etal_2017} or unbounded.

\subparagraph{Alt$(iii)$}

This setting is a cross-hatch of Alts ($i$) and ($ii$). Here we have smaller 
$\theta _{0,i}$ $=$ $.00015$ than Alt($i$) for $i$ $=$ $1,...,10$, for more
components than Alt($i$) but fewer than Alt($ii$). As in previous cases $p$%
-max test power is (nearly) unity when $k_{\delta }$ $=$ $0$, or when $%
k_{\delta }$ $=$ $10$ and $k_{\theta ,n}$ $>>$ $n$. In the presence of
nuisance parameters and $n$ $=$ $100$, power is noticeably smaller: $p$-max
test power is $(.039,.174,.342)$ and $(.058,.253,.450)$ when $k_{\theta ,n}$ 
$=$ $200$ and $482$ respectively. Those values are near or exactly unity
once $n$ $=$ $250$.

The $dbl$-max tests, by comparison, yield trivial power with $n$ $=$ $100$,
under either covariate case, and for any $k_{\delta }$ $\geq $ $0$. The
signal appears to be too weak to detect when $n$ $=$ $100$, in particular
given the shrinkage nature of DBL.

\section{Conclusion\label{sec:conclude}}

We present a class of tests for a high dimensional parameter in a regression
setting with nuisance parameters. We focus on a linear model and least
squares estimation, but the method and theory presented here are
generalizable to a broad class of models and estimators. We show how the
hypotheses can be identified by using many ($k_{\theta ,n}$) potentially
vastly lower dimension parameterizations depending on the existence of
nuisance parameters, and we allow for up to $\ln (k_{\theta ,n})$ $=$ $%
o(n^{1/7})$ low dimension models, depending on covariate assumptions. We use
a test statistic that is the maximum in absolute value of the weighted key
estimated parameters across the many low dimension settings. The lower
dimension helps improve estimation accuracy and allows us to sidestep
sparsity requirements, while a max-statistic both alleviates the need for a
multivariate normalization used in Wald and score statistics, and hones in
the most informative (weighted) model component. Thus, we avoid inverting a
potentially large dimension variance estimator that may be a poor proxy for
the true sampling dispersion. Indeed, our asymptotic theory sidesteps
traditional extreme value theory arguments, while we instead focus on an
approximate p-value computed by parametric wild bootstrap.

In simulation experiments our max-tests generally have good or sharp size.
They generally dominate bootstrapped debiased Lasso based max-tests in terms
of power, and bootstrapped (and normalized) Wald test in terms of both size
and power in low dimension cases.

Future work may ($i$) lean toward allowing for a broader class of
semi-nonparametric models, including models with a high dimensional
parameter $\theta _{0}$ \textit{and} infinite dimensional unknown function $%
h $ and/or high dimensional nuisance term $\delta _{0}$; ($ii$) verify key
asymptotic theory arguments in a general time series setting; and ($iii$)
expand McKeague and Qian's (\citeyear{McKeague_Qian_2015}) significant
predictor method to a high dimensional setting.

\setcounter{equation}{0} \renewcommand{\theequation}{{\thesection}.%
\arabic{equation}} \appendix

\section{Appendix\label{app:a}}

\setstretch{1.2}Throughout $O_{p}(1)$ and $o_{p}(1)$ are not functions of
model counter $i$. $\{k_{\theta ,n},k_{n}\}$ are monotonically increasing
sequences of positive integers. Define parsimonious parameter spaces $%
\mathcal{B}_{(i)}$ $\equiv $ $\mathcal{D}$ $\times $ $\Theta _{(i)}$ where $%
\mathcal{D}$ $\subset $ $\mathbb{R}^{k_{\delta }}$ and $\Theta _{(i)}$ $%
\subset $ $\mathbb{R}$ are compact subsets, $\delta _{0}$ is an interior
point of $\mathcal{D}$, and $0$ and $\theta _{0,i}$ are interior points of $%
\Theta _{(i)}$. Define compact $\mathcal{B}$ $\equiv $ $\mathcal{D}$ $\times 
$ $\Theta $ where $\Theta $ $\equiv $ $\times _{i=1}^{k_{\theta ,n}}\Theta
_{i}$. Recall $\widehat{\mathcal{H}}_{(i)}$ $\equiv $ $1/n%
\sum_{t=1}^{n}x_{(i),t}x_{(i),t}^{\prime }$ and $\mathcal{H}_{(i)}\equiv
E[x_{(i),t}x_{(i),t}^{\prime }]$. Let \underline{$\lambda $}$_{(i),n}$ and 
\underline{$\lambda $}$_{(i)}$\ denote the minimum eigenvalues of $\widehat{%
\mathcal{H}}_{(i)}$ and $\mathcal{H}_{(i)}$

\subsection{Assumptions\label{app:assum}}

\begin{assumption}
\label{assum:smooth} \ \ \medskip \newline
$a.$ $(\epsilon _{t},x_{t})$ are iid over $t$; $E[\epsilon _{t}]$ $=$ $0$; 
\b{c} $\leq $ $E[\epsilon _{t}^{2}]$ $\leq $ $\bar{c}$ and \b{c} $\leq $ $%
E[x_{j,t}^{2}]$ $\leq $ $\bar{c}$\ for all $j$ and\textit{\ some }$\text{\b{c%
},}\bar{c}$ $\in $ $(0,\infty )$ that may differ for different variables; $%
E[\epsilon _{t}^{4}]$ $<$ $\infty $ and $E[x_{j,t}^{4}]$ $<$ $\infty $\ for
all $j$; $P(E[\epsilon _{t}^{2}|x_{t}]$ $=$ $\sigma ^{2})$ $=$ $1$ for
finite $\sigma ^{2}$ $>$ $0$; and $\lim \sup_{n\rightarrow \infty }|\theta
_{0}|$ $<$ $\infty $.$\medskip $\newline
$b.$ $\beta _{0}$ uniquely minimizes $E[(y_{t}$ $-$ $\beta ^{\prime
}x_{t})^{2}]$ on $\mathcal{B}$; $E[(y_{t}$ $-$ $\beta _{(i)}^{\ast \prime
}x_{(i),t})x_{(i),t}]$ $=$ $\boldsymbol{0}_{k_{\delta }+1}$ for all $i$ and
unique $\beta _{(i)}^{\ast }$ in the interior of $\mathcal{B}_{(i)}$.$%
\medskip $\newline
$c.$ $\liminf\limits_{n\rightarrow \infty }\min\limits_{1\leq i\leq
k_{\theta ,n}}\inf\limits_{\lambda ^{\prime }\lambda =1}E[(\lambda ^{\prime
}x_{(i),t})^{2}]$ $>$ $0$ and $\liminf\limits_{n\rightarrow \infty
}\min\limits_{1\leq i\leq k_{\theta ,n}}\{\underline{\lambda }_{(i)}\}$ $>$ $%
0.\medskip $\newline
$d.$ $\liminf\limits_{n\rightarrow \infty }\inf\limits_{\lambda ^{\prime
}\lambda =1}\min\limits_{i\in \mathbb{N}}\{\frac{1}{n}\sum_{t=1}^{n}(\lambda 
\mathcal{H}_{(i)}^{-1\prime }x_{(i),t})^{2}\}$ $>$ $0$ $a.s.$; $%
\liminf\limits_{n\rightarrow \infty }\min\limits_{i\in \mathbb{N}}\{%
\underline{\lambda }_{(i),n}\}$ $>$ $0$ $a.s$.
\end{assumption}

\begin{remark}
\normalfont($a$)-($d$) are standard for linear iid regression models,
augmented to allow for high dimension. The existence of higher moments is
not unusual for a high dimensional maximum limit theory %
\citep[e.g.][]{BuhlmannVanDeGeer2011,Chernozhukov_etal2013,Dezeure_etal_2017,ZhangWu2017}%
. Indeed, in order to handle the Gaussian approximation error $\max_{1\leq
i\leq k_{\theta ,n}}|[\boldsymbol{0}_{k_{\delta }}^{\prime },1]\mathcal{\hat{%
R}}_{i}|$ in (\ref{expand}), an upper bound on feasible $k_{\theta ,n}$ is
generally linked to higher moments. See the discussion following Lemma \ref%
{lm:expansion}.
\end{remark}

\begin{remark}
\normalfont By the constructions of the $\mathcal{L}_{2}$-minimizers $\beta
_{0}$ and $\beta _{(i)}^{\ast }$, (\ref{model}) and the parsimonious
versions may be misspecified in the senses of \cite{Sawa1978} and \cite%
{White1982}.
\end{remark}

\begin{remark}
\normalfont We use the eigenvalue bounds $\lim \inf_{n\rightarrow \infty
}\min_{1\leq i\leq k_{\theta ,n}}\{\underline{\lambda }_{(i),n}\}>$ $0$ $%
a.s. $ to ensure $\max_{1\leq i\leq k_{\theta ,n}}||\widehat{\mathcal{H}}%
_{(i)}^{-1}||$ $=$ $\max_{1\leq i\leq k_{\theta ,n}}\{\underline{\lambda }%
_{(i),n}^{-1}\}$ $=$ $(\min_{1\leq i\leq k_{\theta ,n}}\{\underline{\lambda }%
_{(i),n}^{-1}\})^{-1}$ $=$ $O_{p}(1)$ is well defined for each $n$ and any $%
\{k_{\theta ,n}\}$. The condition $|\theta _{0}|$ $<$ $\infty $ uniformly in 
$n$\ is only used in Lemma \ref{lm:beta_hat_rate} in order to bound the
parsimonious model error moment $E[\max_{1\leq i\leq k_{\theta
,n}}\left\vert v_{(i),t}\right\vert ^{4}]$ under \emph{either} hypothesis.
Under the null, of course, $E[\max_{1\leq i\leq k_{\theta ,n}}\left\vert
v_{(i),t}\right\vert ^{4}]$ $=$ $E|\epsilon _{t}|^{4}$ $<$ $\infty $.
\end{remark}

\subsection{Proofs of Main Results\label{app:main_proofs}}

Write compactly $|\cdot |_{\infty }$ $\equiv $ $\max_{1\leq i\leq k_{\theta
,n}}|\cdot |$ and $||\cdot ||_{\infty }$ $\equiv $ $\max_{1\leq i\leq
k_{\theta ,n}}||\cdot ||$.\medskip \newline
\textbf{Proof of Theorem \ref{th:ident}.} By supposition $E[(y_{t}$ $-$ $%
\beta ^{\prime }x_{t})x_{i,t}]$ $=$ $0$ $\forall i$ \emph{if and only if} $%
\beta =\beta _{0}=\left[ \delta _{0}^{\prime },\theta _{0}^{\prime }\right]
^{\prime }$. Therefore, if $\theta _{0}$ $=$ $\boldsymbol{0}_{k_{\theta ,n}}$
then $E[(y_{t}$ $-$ $\delta _{0}^{\prime }x_{\delta ,t})x_{i,t}]=0$ $\forall
i$ hence $E[(y_{t}$ $-$ $\delta _{0}^{\prime }x_{\delta ,t})x_{(i),t}]=%
\boldsymbol{0}_{k_{\delta }+1}$ for each $i$ $=$ $1,...,k_{\theta ,n}$. This
instantly yields $\theta _{i}^{\ast }$ $=$ $0$ for each $i$ and $\delta
^{\ast }$ $=$ $\delta _{0}$\ by construction and uniqueness of $\beta
_{(i)}^{\ast }$ $=$ $[\delta ^{\ast \prime },\theta _{i}^{\ast }]^{\prime }$.

Conversely, if $\theta ^{\ast }$ $=$ $\boldsymbol{0}_{k_{\theta ,n}}$ then
from the parsimonious risk it follows $E[(y_{t}$ $-$ $\delta _{(i)}^{\ast
\prime }x_{\delta ,t})x_{(i),t}]=\boldsymbol{0}_{k_{\delta }+1}$ for each $i$
by identification of $\beta _{(i)}^{\ast }$ $=$ $[\delta ^{\ast \prime
},\theta _{i}^{\ast }]^{\prime }$. But this implies the sub-gradient $%
E[(y_{t}$ $-$ $\delta _{(i)}^{\ast \prime }x_{\delta ,t})x_{\delta ,t}]=%
\boldsymbol{0}_{k_{\delta }+1}$ for each $i$, hence the $\delta _{(i)}^{\ast
\prime }s$ must be identical: there exists a unique $\delta ^{\ast }$ $\in $ 
$\mathcal{D}$ such that $\delta _{(i)}^{\ast }$ $=$ $\delta ^{\ast }$ for
each $i$. But $E[(y_{t}$ $-$ $\delta ^{\ast \prime }x_{\delta ,t})x_{(i),t}]$
$=$ $\boldsymbol{0}_{k_{\delta }+1}$ $\forall i$ yields $E[(y_{t}$ $-$ $%
\delta ^{\ast \prime }x_{\delta ,t})x_{j,t}]=0$ $\forall j$ by the
definition of $x_{(i),t}$. Therefore $[\delta _{0},\theta _{0}]$ $=$ $%
[\delta ^{\ast },\boldsymbol{0}_{k_{\theta ,n}}]$ by identification of $%
[\delta _{0},\theta _{0}]$. $\mathcal{QED}$.\medskip

\noindent \textbf{Proof of Theorem \ref{thm:max_dist}.} Define\footnote{%
We suppress $\lambda $, $\lambda ^{\prime }\lambda $ $=$ $1$, when confusion
is avoided.} $s_{(i),t}$ $=$ $s_{(i),t}(\lambda )$ $\equiv $ $\lambda
^{\prime }(E[x_{(i),t}x_{(i),t}^{\prime }])^{-1}x_{(i),t}$; $\mathcal{M}_{r}$
$\equiv $ $||\epsilon _{t}||_{r}$ $\times $ $\max_{1\leq i\leq k_{\theta
,n}}\left\Vert s_{(i),t}\right\Vert _{r}$ for $r$ $\in $ $(2,4]$, where $%
\lim \sup_{n\rightarrow \infty }\mathcal{M}_{r}$ $<$ $\infty $ under
Assumption \ref{assum:smooth}; $\zeta _{(i),t}$ $\equiv $ $\epsilon
_{t}s_{(i),t}/(||\epsilon _{t}||_{2}||s_{(i),t}||_{2})$; define random
variables $z_{(i),t}$ $\sim $ $N(0,1)$ independent over $t$; and for any $%
\gamma $ $\in $ $(0,1)$ define $u(\gamma )$ $\equiv $ $u_{\epsilon s}(\gamma
)$ $\vee $ $u_{z}(\gamma )$ where 
\begin{eqnarray*}
&&u_{\epsilon s}(\gamma )\equiv \inf \left\{ u:P\left( \max_{1\leq i\leq
k_{\theta ,n},1\leq t\leq n}\left\vert \zeta _{(i),t}\right\vert >u\right)
\leq \gamma \right\} \\
&&u_{z}(\gamma )\equiv \inf \left\{ u:P\left( \max_{1\leq i\leq k_{\theta
,n},1\leq t\leq n}\left\vert z_{(i),t}\right\vert >u\right) \leq \gamma
\right\} .
\end{eqnarray*}

Write $\mathcal{K}_{3,4}$ $\equiv $ $\mathcal{M}_{3}^{3/4}$ $\vee $ $%
\mathcal{M}_{4}^{1/2}$. Under Assumption \ref{assum:smooth}.a,c and $H_{0}$,
Theorem 2.2 in \cite{Chernozhukov_etal2013} applies to $\epsilon
_{t}s_{(i),t}$. Hence, for every $\gamma $ $\in $ $(0,1)$, and some finite
constant $\mathcal{C}$ $>$ $0$ that depends only on lower and upper bounds
on $E[\epsilon _{t}^{2}]$ and $E[s_{j,t}^{2}]$, 
\begin{equation}
\rho _{n}(\lambda )\leq \mathcal{C}\left[ \frac{1}{n^{1/8}}\mathcal{K}%
_{3,4}\left( \ln \left( \frac{k_{\theta ,n}n}{\gamma }\right) \right) ^{7/8}+%
\frac{1}{n^{1/2}}\left( \ln \left( \frac{k_{\theta ,n}n}{\gamma }\right)
\right) ^{3/2}u(\gamma )+\gamma \right] .  \label{rho_n}
\end{equation}

We need to bound $u(\gamma )$. Notice 
\begin{equation*}
\sup_{\lambda ^{\prime }\lambda =1}E\left[ \max_{1\leq i\leq k_{\theta
,n}}s_{(i),t}^{4}(\lambda )\right] \leq \left\vert \left(
E[x_{(i),t}x_{(i),t}^{\prime }]\right) ^{-1}\right\vert _{\infty
}^{4}E\left\vert x_{(i),t}\right\vert _{\infty }^{4}=\left\vert \mathcal{H}%
_{(i)}^{-1}\right\vert _{\infty }^{4}\mathcal{M}_{n},
\end{equation*}%
where $\mathcal{M}_{n}$ $\equiv $ $E|x_{(i),t}|_{\infty }^{4}$, and $|%
\mathcal{H}_{(i)}^{-1}|_{\infty }$ $\in $ $(0,\infty )$ by Assumption \ref%
{assum:smooth}.c. Assumption \ref{assum:smooth}.c ensures $\chi $ $\equiv $ $%
\liminf_{n\rightarrow \infty }\min_{1\leq i\leq k_{\theta ,n}}\inf_{\lambda
^{\prime }\lambda =1}||s_{(i),t}(\lambda )||_{2}$ $>$ $0$. Hence, by the
Cauchy-Schwartz inequality:%
\begin{eqnarray*}
\left( E\left[ \max_{1\leq i\leq k_{\theta ,n}}\zeta _{(i),t}^{2}\right]
\right) ^{1/2} &\leq &\frac{\left\Vert \epsilon _{t}\right\Vert _{4}}{%
\left\Vert \epsilon _{t}\right\Vert _{2}}\left( E\left[ \max_{1\leq i\leq
k_{\theta ,n}}\frac{s_{(i),t}^{4}}{\left\Vert s_{(i),t}\right\Vert _{2}^{4}}%
\right] \right) ^{1/4} \\
&\leq &\frac{\left\Vert \epsilon _{t}\right\Vert _{4}}{\chi \left\Vert
\epsilon _{t}\right\Vert _{2}}\left\Vert \mathcal{H}_{(i)}^{-1}\right\Vert
_{\infty }^{1/4}\times \mathcal{M}_{n}^{1/4}\equiv \mathcal{K\times M}%
_{n}^{1/4},
\end{eqnarray*}%
where $\mathcal{K}$\ is implicit. Hence, $u_{\epsilon s}(\gamma )$ $=$ $%
\mathcal{KM}_{n}^{1/4}/\gamma $ in view of Markov and Liapunov inequalities:%
\begin{eqnarray*}
P\left( \max_{1\leq i\leq k_{\theta ,n},1\leq t\leq n}\left\vert \zeta
_{(i),t}\right\vert >u_{\epsilon s}\right) &\leq &\frac{1}{u_{\epsilon s}}E%
\left[ \max_{1\leq i\leq k_{\theta ,n},1\leq t\leq n}\left\vert \zeta
_{(i),t}\right\vert \right] \\
&\leq &\frac{1}{u_{\epsilon s}}\left( E\left[ \max_{1\leq i\leq k_{\theta
,n},1\leq t\leq n}\zeta _{(i),t}^{2}\right] \right) ^{1/2}\leq \frac{1}{%
u_{\epsilon s}}\mathcal{KM}_{n}^{1/4}.\text{ }
\end{eqnarray*}%
Second, use a log-exp argument with Jensen's inequality, and Gaussianicity,
to yield 
\begin{eqnarray*}
E\left[ \max_{1\leq i\leq k_{\theta ,n},1\leq t\leq n}\left\vert
z_{(i),t}\right\vert \right] &\leq &\frac{1}{\lambda }\ln \left( E\left[
\exp \left( \lambda \max_{1\leq i\leq k_{\theta ,n},1\leq t\leq n}\left\vert
z_{(i),t}\right\vert \right) \right] \right) \text{ for any }\lambda >0: \\
&\leq &\frac{1}{\lambda }\ln \left( nk_{\theta ,n}E\left[ \exp \left\{
\lambda \left\vert z_{(i),t}\right\vert \right\} \right] \right) \\
&=&\frac{1}{\lambda }\ln \left( nk_{\theta ,n}\right) +\lambda =2\sqrt{\ln
\left( nk_{\theta ,n}\right) }.
\end{eqnarray*}%
The final equality uses the minimizer $\sqrt{\ln \left( nk_{\theta
,n}\right) }$ of $\lambda ^{-1}\ln \left( nk_{\theta ,n}\right) $ $+$ $%
\lambda $. Hence $u_{z}(\gamma )$ $=$ $(2/\gamma )\sqrt{\ln \left(
nk_{\theta ,n}\right) }$. Combine the above to yield:%
\begin{equation*}
u(\gamma )=\frac{1}{\gamma }\left\{ \mathcal{KM}_{n}^{1/4}\vee 2\sqrt{\ln
\left( nk_{\theta ,n}\right) }\right\} .
\end{equation*}%
Now use cases $(i)$-$(iv)$ for $\mathcal{M}_{n}$ derived in the proof of
Lemma \ref{lm:expansion} to deduce bounds for $u(\gamma )$, for some
universal $\mathcal{C}$ $>$ $0$: $(i)$ $\mathcal{M}_{n}$ $\leq $ $K$ hence $%
u(\gamma )$ $\leq $ $\mathcal{C}\frac{1}{\gamma }\sqrt{\ln \left( nk_{\theta
,n}\right) }$; $(ii)$ $\mathcal{M}_{n}$ $\leq $ $K\sqrt{\ln (k_{\theta ,n})}$
hence $u(\gamma )$ $\leq $ $\frac{1}{\gamma }\mathcal{C}\sqrt{\ln \left(
nk_{\theta ,n}\right) }$; $(iii)$ $\mathcal{M}_{n}$ $\leq $ $K\ln (k_{\theta
,n})$ hence again $u(\gamma )$ $=$ $\frac{1}{\gamma }\mathcal{C}\sqrt{\ln
\left( nk_{\theta ,n}\right) }$; and $(iv)$ $\mathcal{M}_{n}\leq Kk_{\theta
,n}^{4/p}$ yields $u(\gamma )$ $\leq $ $\frac{1}{\gamma }\{\mathcal{K}%
k_{\theta ,n}^{1/p}$ $\vee $ $2\sqrt{\ln \left( nk_{\theta ,n}\right) }\}$ $%
\leq $ $\mathcal{C}\frac{1}{\gamma }k_{\theta ,n}^{1/p}\sqrt{\ln \left(
n\right) }$.

Next, return to (\ref{rho_n}). Let $\{g(n)\}$ be a positive monotonic
sequence, $g(n)$ $\rightarrow $ $\infty $. Then $\rho _{n}(\lambda )$ $=$ $%
O(1/g(n))$ $=$ $o(1)$ for any $\gamma $ $=$ $O(1/g(n))$ whenever $%
n^{-1/8}\left( \ln \left( k_{\theta ,n}n/\gamma \right) \right) ^{7/8}$ $%
\rightarrow $ $0$ and $n^{-1/2}\left( \ln \left( k_{\theta ,n}n/\gamma
\right) \right) ^{3/2}u(\gamma )$ $\rightarrow $ $0$. Using $\gamma $ $=$ $%
O(1/g(n))$, the first condition reduces to 
\begin{equation}
\ln (k_{\theta ,n})+\ln \left( ng\left( n\right) \right)
=o(n^{1/7})\Longrightarrow \ln (k_{\theta ,n})=o(n^{1/7})\text{ if }\ln
(g(n))=o(n^{1/7}).  \label{kn1}
\end{equation}%
The second condition reduces to $n^{-1/2}[\ln \left( nk_{\theta ,n}\right) $ 
$+$ $\ln (g\left( n\right) ]^{3/2}u(\gamma )$ $=$ $o\left( 1\right) $.
Hence, for any slowly varying $g\left( n\right) $ $\rightarrow $ $\infty $,
by the above case-specific bounds for $u(\gamma )$:%
\begin{eqnarray}
&&\text{ }(i)\text{-}(iii)\ln \left( nk_{\theta ,n}\right) =o\left(
n^{1/4}/g(n)^{1/2}\right) \Longrightarrow \ln (k_{\theta ,n})=o\left( \frac{%
n^{1/4}}{g(n)^{1/2}}\right)  \label{kn2} \\
&&\text{ }(iv)\text{ }\ln \left( n\right) ^{3/2}\left\{ 1+\frac{\ln \left(
k_{\theta ,n}g\left( n\right) \right) }{\ln \left( n\right) }\right\}
^{3/2}k_{\theta ,n}^{1/p}=o\left( \frac{n^{1/2}}{g(n)\sqrt{\ln \left(
n\right) }}\right) \Longrightarrow k_{\theta ,n}=o\left( \frac{n^{p/2}}{%
\left\{ g(n)\ln ^{2}\left( n\right) \right\} ^{p}}\right)  \notag
\end{eqnarray}%
Now take the minimum $k_{\theta ,n}$ from (\ref{kn1}) and (\ref{kn2}) to
yield $\rho _{n}(\lambda )$ $=$ $O(1/g(n))$ $\rightarrow $ $0$ for slowly
varying $g\left( n\right) $ $\rightarrow $ $\infty $ provided by case: $(i)$-%
$(iii)$ $\ln (k_{\theta ,n})=o(n^{1/7})$, and $(iv)$ $k_{\theta ,n}$ $=$ $%
o(n^{p/2}/[g(n)\ln ^{2}\left( n\right) ]^{p})$. $\mathcal{QED}$.\medskip

The proof of Theorem \ref{th:max_theta_hat}, and multiple lemmas in SM,
exploit a general high dimensional probability bound applied to $\widehat{%
\mathcal{H}}_{(i)}$ $-$ $\mathcal{H}_{(i)}$, $\sqrt{n}\widehat{\mathcal{G}}%
_{(i)}$ and $\sqrt{n}(\hat{\beta}_{(i)}$ $-$ $\beta _{(i)}^{\ast })$. Recall 
$\widehat{\mathcal{G}}_{(i)}$ $\equiv $ $-1/n%
\sum_{t=1}^{n}v_{(i),t}x_{(i),t} $, $\widehat{\mathcal{H}}_{(i)}$ $\equiv $ $%
1/n\sum_{t=1}^{n}x_{(i),t}x_{(i),t}^{\prime }$, $\mathcal{H}_{(i)}$ $\equiv $
$E[x_{(i),t}x_{(i),t}^{\prime }]$, and $\mathcal{M}_{n}$ $\equiv $ $%
E|x_{(i),t}|_{\infty }^{4}$.

\begin{lemma}
\label{lm:beta_hat_rate}Let Assumption \ref{assum:smooth} hold, and let $%
\{k_{\theta ,n}\}$ be any monotonically increasing sequence of integers.
Then $(a)$ $||\sqrt{n}\widehat{\mathcal{G}}_{(i)}||_{\infty }$ $=$ $O_{p}(%
\sqrt{\ln (k_{\theta ,n})\mathcal{M}_{n}})$; $(b)$ $||\widehat{\mathcal{H}}%
_{(i)}^{-1}$ $-$ $\mathcal{H}_{(i)}^{-1}||_{\infty }$ $=$ $O_{p}(\sqrt{\ln
(k_{\theta ,n})\mathcal{M}_{n}/n})$; and $(c)$ $||\sqrt{n}(\hat{\beta}_{(i)}$
$-$ $\beta _{(i)}^{\ast })||_{\infty }$ $=O_{p}(\sqrt{\ln (k_{\theta ,n})%
\mathcal{M}_{n}})$ provided $\ln (k_{\theta ,n})$ $=$ $o(\sqrt{n}/\mathcal{M}%
_{n})$.
\end{lemma}

\begin{remark}
\normalfont$\ln (k_{\theta ,n})$ $=$ $o(\sqrt{n}/\mathcal{M}_{n})$ holds
under each covariate case $(i)$ $-$ $(iv)$: see the proof of Lemma \ref%
{lm:expansion} in SM. Hence $\ln (k_{\theta ,n})\mathcal{M}_{n}/n$ $%
\rightarrow $ $0$, yielding, e.g., $||\widehat{\mathcal{H}}_{(i)}^{-1}$ $-$ $%
\mathcal{H}_{(i)}^{-1}||_{\infty }$ $\overset{p}{\rightarrow }$ $0$ and $||%
\hat{\beta}_{(i)}$ $-$ $\beta _{(i)}^{\ast }||_{\infty }$ $\overset{p}{%
\rightarrow }$ $0$.
\end{remark}

\noindent \textbf{Proof of Theorem \ref{th:max_theta_hat}.} For ($a$),
convergence in distribution (\ref{PP0}) follows from asymptotic expansion
Lemma \ref{lm:expansion}, and convergence in Kolmogorov distance (\ref%
{Gauss_apprx}), cf. Theorem \ref{thm:max_dist} with $\lambda $ $=$ $[%
\boldsymbol{0}_{k_{\delta }}^{\prime },1]^{\prime }$. Then $|\sqrt{n}%
\mathcal{W}_{n,i}\hat{\theta}_{i}|_{\infty }$ $\overset{d}{\rightarrow }$ $%
\max_{i\in \mathbb{N}}|\mathcal{W}_{i}\boldsymbol{Z}_{(i)}([\boldsymbol{0}%
_{k_{\delta }}^{\prime },1])|$ follows immediately. Claim ($b$) follows from
Lemma \ref{lm:beta_hat_rate}.c consistency of $\hat{\theta}_{i}$, and $%
\mathcal{W}_{n,i}$ $\overset{p}{\rightarrow }$ $\mathcal{W}_{i}$ with $%
\mathcal{W}_{i}$ $\in $ $(0,\infty )$. $\mathcal{QED}.\medskip $\newline
\textbf{Proof of Theorem \ref{th:p_value_boot}}.\medskip \newline
\textbf{Claim (a).}\qquad We have $|\sqrt{n}\mathcal{W}_{n,i}\hat{\theta}%
_{i}|_{\infty }$ $\overset{d}{\rightarrow }$ $\max_{i\in \mathbb{N}}|%
\mathcal{W}_{i}\boldsymbol{Z}_{(i)}([\boldsymbol{0}_{k_{\delta }}^{\prime
},1])|$ from Theorem \ref{th:max_theta_hat}.a. Theorem \ref{thm:boot_theta}
yields $\max_{1\leq i\leq k_{\theta ,n}}\mathcal{\tilde{T}}_{n}$ $%
\Rightarrow ^{p}$ $\max_{i\in \mathbb{N}}|\mathcal{W}_{i}\boldsymbol{\tilde{Z%
}}_{(i)}([\boldsymbol{0}_{k_{\delta }}^{\prime },1])|$, where $\{\boldsymbol{%
\tilde{Z}}_{(i)}(\lambda )\}_{i\in \mathbb{N}}$ is an independent copy of $\{%
\boldsymbol{Z}_{(i)}(\lambda )\}_{i\in \mathbb{N}}$, $\boldsymbol{Z}%
_{(i)}(\lambda )$ $\sim $ $N(0,E[\epsilon _{t}^{2}]\lambda ^{\prime }%
\mathcal{H}_{(i)}^{-1}\lambda )$, independent of the asymptotic draw $%
\{x_{t},y_{t}\}_{t=1}^{\infty }$. The proof now follows from arguments in %
\citet[p. 427]{Hansen1996}.\medskip \newline
\textbf{Claim (b).}\qquad Let $H_{1}$ hold, suppose $\theta _{0,j}$ $\neq $ $%
0$ for some $j$, define $\mathfrak{S}_{n}$ $\equiv $ $\{x_{t},y_{t}%
\}_{t=1}^{n}$, and write $\mathcal{P}_{\mathfrak{S}_{n}}(\mathcal{A})$ $%
\equiv $ $\mathcal{P}(\mathcal{A}|\mathfrak{S}_{n})$. In view of $%
\max_{1\leq i\leq k_{\theta ,n}}\mathcal{\tilde{T}}_{n}$ $\Rightarrow ^{p}$ $%
\max_{i\in \mathbb{N}}|\mathcal{W}_{i}\boldsymbol{\tilde{Z}}_{(i)}([%
\boldsymbol{0}_{k_{\delta }}^{\prime },1])|$ we have 
\citep[eq.
(3.4)]{GineZinn1990} 
\begin{equation*}
\sup_{z>0}\left\vert \mathcal{P}_{\mathfrak{S}_{n}}\left( \mathcal{\tilde{T}}%
_{n}\leq z\right) -P\left( \left\vert \mathcal{W}_{i}\boldsymbol{\mathring{Z}%
}_{(i)}([\boldsymbol{0}_{k_{\delta }}^{\prime },1])\right\vert _{\infty
}\leq z\right) \right\vert \overset{p}{\rightarrow }0,
\end{equation*}%
where $\{\boldsymbol{\mathring{Z}}_{(i)}(\cdot )\}_{i\in \mathbb{N}}$ is an
independent copy of $\{\boldsymbol{Z}_{(i)}(\cdot )\}_{i\in \mathbb{N}}$.
Moreover, $\tilde{p}_{n,\mathcal{R}_{n}}$ $=$ $\mathcal{P}_{\mathfrak{S}%
_{n}}(\mathcal{\tilde{T}}_{n,1}$ $>$ $\mathcal{T}_{n})$ $+$ $o_{p}(1)$ by
the Glivenko-Cantelli theorem and independence across bootstrap samples.
Further, $\mathcal{T}_{n}$ $\geq $ $|\sqrt{n}\mathcal{W}_{n,j}(\hat{\theta}%
_{j}$ $-$ $\theta _{0,j})+\sqrt{n}\mathcal{W}_{n,j}\theta _{0,j}|$ $=$ $%
\sqrt{n}\mathcal{W}_{j}\left\vert \theta _{0,j}\right\vert (|\hat{\theta}%
_{j}/\theta _{0,j}|$ $+$ $o_{p}(1))$ $\overset{p}{\rightarrow }$ $\infty $
by Theorem \ref{th:max_theta_hat}.b and $|\mathcal{W}_{n,i}$ $-$ $\mathcal{W}%
_{i}|_{\infty }$ $=$ $o_{p}(1)$. Therefore%
\begin{equation*}
\tilde{p}_{n,\mathcal{R}_{n}}\leq P\left( \left\vert \mathcal{W}_{i}%
\boldsymbol{\mathring{Z}}_{(i)}([\boldsymbol{0}_{k_{\delta }}^{\prime
},1])\right\vert _{\infty }>\sqrt{n}\mathcal{W}_{j}\left\vert \theta
_{0,j}\right\vert \left( \left\vert \frac{\hat{\theta}_{j}}{\theta _{0,j}}%
\right\vert +o_{p}(1)\right) \right) +o_{p}(1)\rightarrow 0
\end{equation*}%
as long as $|\mathcal{W}_{i}\boldsymbol{\mathring{Z}}_{(i)}([\boldsymbol{0}%
_{k_{\delta }}^{\prime },1])|_{\infty }$ $=$ $o_{p}(\sqrt{n})$. We show
below that, given Gaussianicity of $\boldsymbol{\mathring{Z}}_{(i)}([%
\boldsymbol{0}_{k_{\delta }}^{\prime },1])$: 
\begin{equation}
\left\vert \mathcal{W}_{i}\boldsymbol{\mathring{Z}}_{(i)}([\boldsymbol{0}%
_{k_{\delta }}^{\prime },1])\right\vert _{\infty }=O_{p}\left( \sqrt{\ln
\left( k_{\theta ,n}\right) }\right)  \label{PWZ}
\end{equation}%
hence $\tilde{p}_{n,\mathcal{R}_{n}}$ $\overset{p}{\rightarrow }$ $0$ under $%
H_{1}$ and therefore $P(\tilde{p}_{n,\mathcal{R}_{n}}$ $<$ $\alpha )$ $%
\rightarrow $ $1$ as claimed.

We now prove (\ref{PWZ}). Put $\zeta $ $=$ $\sqrt{2\mathcal{W}E[\epsilon
_{t}^{2}]\mathcal{H}^{-1,\theta ,\theta }}$ $>$ $0$ where $\mathcal{W}$ $%
\equiv $ $\max_{i\in \mathbb{N}}\{\mathcal{W}_{i}\}$ $<$ $\infty $, and $%
\mathcal{H}^{-1,\theta ,\theta }$ $\equiv $ $\max\nolimits_{1\leq i\leq
k_{\theta ,n}}[\boldsymbol{0}_{k_{\delta }}^{\prime },1]\mathcal{H}%
_{(i)}^{-1}[\boldsymbol{0}_{k_{\delta }}^{\prime },1]^{\prime }$ $\in $ $%
(0,\infty )$ under positive definiteness Assumption \ref{assum:smooth}.c.
Then Chernoff's Gaussian concentration bound $P(|Y|$ $>$ $c)$ $\leq $ $%
e^{-c^{2}/(2\sigma ^{2})}$ for $Y$ $\sim $ $N(0,\sigma ^{2})$, and Boole's
inequality yield for any $\xi $ $>$ $0$:%
\begin{eqnarray}
&&P\left( \left\vert \mathcal{W}_{i}\boldsymbol{\mathring{Z}}_{(i)}([%
\boldsymbol{0}_{k_{\delta }}^{\prime },1])\right\vert _{\infty }>\zeta \sqrt{%
\ln \left( k_{\theta ,n}\right) }+\xi \right)  \label{PWZ1} \\
&&\text{ \ \ \ \ \ \ \ \ \ \ \ \ \ \ \ }\leq \sum_{i=1}^{k_{\theta
,n}}P\left( \left\vert \mathcal{W}_{i}\boldsymbol{\mathring{Z}}_{(i)}([%
\boldsymbol{0}_{k_{\delta }}^{\prime },1])\right\vert >\zeta \sqrt{\ln
\left( k_{\theta ,n}\right) }+\xi \right)  \notag \\
&&\text{ \ \ \ \ \ \ \ \ \ \ \ \ \ \ \ }\leq \sum_{i=1}^{k_{\theta ,n}}\exp
\left\{ -\frac{\left( \zeta \sqrt{\ln \left( k_{\theta ,n}\right) }+\xi
\right) ^{2}}{2\mathcal{W}E[\epsilon _{t}^{2}]\mathcal{H}^{-1,\theta ,\theta
}}\right\}  \notag \\
&&\text{ \ \ \ \ \ \ \ \ \ \ \ \ \ \ \ }=\frac{1}{k_{\theta ,n}}%
\sum_{i=1}^{k_{\theta ,n}}\exp \left\{ -\frac{2\xi \sqrt{\ln \left(
k_{\theta ,n}\right) }}{\zeta }\right\} \exp \left\{ -\frac{\xi ^{2}}{\zeta
^{2}}\right\} \leq \exp \left\{ -\xi ^{2}/\zeta ^{2}\right\} .\text{ }%
\mathcal{QED}.  \notag
\end{eqnarray}%
\textbf{Proof of Theorem \ref{th:local_power}.} Write $\xi _{n}$ $\equiv $ $%
\sqrt{n/\ln \left( k_{\theta ,n}\right) \mathcal{M}_{n}}$. Under $H_{1}^{L}$
by the same arguments used to prove Lemma \ref{lm:expansion} and Theorem \ref%
{th:max_theta_hat}, and in view of Lemma \ref{lm:beta_hat_rate}.c: 
\begin{equation}
\left\vert \left\vert \sqrt{n}\mathcal{W}_{n,i}\left( \hat{\theta}%
_{i}-c_{i}\xi _{n}\right) \right\vert _{\infty }-\left\vert \mathcal{W}_{i}%
\boldsymbol{Z}_{(i)}([\boldsymbol{0}_{k_{\delta }}^{\prime },1])\right\vert
_{\infty }\right\vert \overset{p}{\rightarrow }0,  \label{theta_HL1}
\end{equation}%
where $\boldsymbol{Z}_{(i)}([\boldsymbol{0}_{k_{\delta }}^{\prime },1])$ $%
\sim $ $N(0,E[\epsilon _{t}^{2}]\mathcal{H}^{-1,\theta ,\theta })$ with $%
\mathcal{H}^{-1,\theta ,\theta }$ $\equiv $ $\max\nolimits_{1\leq i\leq
k_{\theta ,n}}[\boldsymbol{0}_{k_{\delta }}^{\prime },1]\mathcal{H}%
_{(i)}^{-1}[\boldsymbol{0}_{k_{\delta }}^{\prime },1]^{\prime }$.

Now define for arbitrary $d$ $=$ $[d_{i}]_{i=1}^{k_{\theta ,n}}$ and sample $%
\mathfrak{S}_{n}$ $=$ $\{x_{t},y_{t}\}_{t=1}^{n}$: 
\begin{equation*}
\mathcal{P}_{n}(d)\equiv P\left( \frac{\left\vert \mathcal{W}_{i}\boldsymbol{%
\mathring{Z}}_{(i)}([\boldsymbol{0}_{k_{\delta }}^{\prime },1])\right\vert
_{\infty }}{\sqrt{\ln \left( k_{\theta ,n}\right) \mathcal{M}_{n}}}%
>\left\vert \xi _{n}\mathcal{W}_{n,i}\left( \hat{\theta}_{i}-c_{i}/\xi
_{n}\right) +d_{i}\right\vert _{\infty }|\mathfrak{S}_{n}\right) .
\end{equation*}%
where $\{\boldsymbol{\mathring{Z}}_{(i)}(\cdot )\}_{i\in \mathbb{N}}$ is an
independent copy of $\{\boldsymbol{Z}_{(i)}(\cdot )\}_{i\in \mathbb{N}}$. By
the Glivenko-Cantelli Theorem with $\mathcal{R}$ $=$ $\mathcal{R}_{n}$ $%
\rightarrow $ $\infty $, the triangle inequality, $|\mathcal{W}_{n,i}$ $-$ $%
\mathcal{W}_{i}|_{\infty }$ $=$ $O_{p}(\sqrt{\ln \left( k_{\theta ,n}\right) 
\mathcal{M}_{n}/n})$, and the stated bounds on $k_{\theta ,n}$ which ensure $%
\ln \left( k_{\theta ,n}\right) $ $=$ $o(\sqrt{n}/\mathcal{M}_{n})$ in all
cases, and arguments in the proof of Theorem \ref{th:p_value_boot}.b, the
bootstrapped p-value approximation satisfies:%
\begin{eqnarray*}
\tilde{p}_{n,\mathcal{R}} &=&P\left( \left\vert \sqrt{n}\mathcal{W}_{n,i}%
\widehat{\tilde{\theta}}_{i}\right\vert _{\infty }>\left\vert \sqrt{n}%
\mathcal{W}_{n,i}\left( \hat{\theta}_{i}-c_{i}\xi _{n}\right) +\mathcal{W}%
_{i}c_{i}\sqrt{\ln \left( k_{\theta ,n}\right) \mathcal{M}_{n}}\right\vert
_{\infty }|\mathfrak{S}_{n}\right) +o_{p}(1) \\
&=&P\left( \frac{\left\vert \mathcal{W}_{i}\boldsymbol{\mathring{Z}}_{(i)}([%
\boldsymbol{0}_{k_{\delta }}^{\prime },1])\right\vert _{\infty }}{\sqrt{\ln
\left( k_{\theta ,n}\right) \mathcal{M}_{n}}}>\left\vert \xi _{n}\mathcal{W}%
_{n,i}\left( \hat{\theta}_{i}-c_{i}/\xi _{n}\right) +\mathcal{W}%
_{i}c_{i}\right\vert _{\infty }|\mathfrak{S}_{n}\right) +o_{p}(1) \\
&=&\mathcal{P}_{n}\left( \mathcal{W}_{j}c_{j}\right) +o_{p}(1).
\end{eqnarray*}%
The proof of Theorem \ref{th:p_value_boot}.a and arguments in 
\citet[p.
427]{Hansen1996} imply $\mathcal{P}_{n}(\boldsymbol{0}_{k_{\theta ,n}})$ is
asymptotically uniformly distributed.

Now, in view of (\ref{PWZ}) and (\ref{theta_HL1}) it follows that $%
\max_{1\leq i\leq k_{\theta ,n}}|\mathcal{W}_{i}\boldsymbol{\mathring{Z}}%
_{(i)}([\boldsymbol{0}_{k_{\delta }}^{\prime },1])|$ $=$ $O_{p}(\sqrt{\ln
\left( k_{\theta ,n}\right) })$ and $\max_{1\leq i\leq k_{\theta ,n}}|\xi
_{n}\mathcal{W}_{n,i}(\hat{\theta}_{i}$ $-$ $c_{i}/\xi _{n})|$ $=$ $O_{p}(1)$%
. Thus $\lim_{n\rightarrow \infty }P(\tilde{p}_{n,\mathcal{R}_{n}}$ $<$ $%
\alpha )$ $\nearrow $ $1$ monotonically as the maximum local drift $%
\max_{i\in \mathbb{N}}|c_{i}|$ $\nearrow $ $\ \infty $ since $\mathcal{W}%
_{i} $ $\in $ $(0,\infty )$ $\forall i$. $\mathcal{QED}$

\setstretch{.5} 
\bibliographystyle{econometrica}
\bibliography{refs_maxtest_many_zeros_LINEAR}

\singlespacing\setstretch{1} \clearpage
\newpage \clearpage

\enlargethispage{2cm}

\begin{sidewaystable}[!t]
\caption{Rejection Frequencies under $H_{0}$ and $H_{1}(i)$} \label%
{table:max_wald_rejH0H11}

\begin{center}
{\small 
\begin{tabular}{l|lll|lll|lll|lll|lll|lll}
\hline
& \multicolumn{9}{|c|}{$H_{0}:\theta _{0}=0$} & \multicolumn{9}{|c}{$%
H_{1}(i):\theta _{0,1}$ $=$ $.0015$ with $[\theta _{0,i}]_{i=2}^{k_{\theta
,n}}$ $=$ $\boldsymbol{0}_{k_{\theta ,n}-1}$} \\ \hline
\multicolumn{1}{c|}{} & \multicolumn{9}{|c|}{$n=100$} & \multicolumn{9}{|c}{$%
n=100$} \\ \hline\hline
& \multicolumn{3}{|c|}{$k_{\theta ,n}=35$} & \multicolumn{3}{|c|}{$k_{\theta
,n}=200$} & \multicolumn{3}{|c|}{$k_{\theta ,n}=482$} & \multicolumn{3}{|c|}{%
$k_{\theta ,n}=35$} & \multicolumn{3}{|c|}{$k_{\theta ,n}=200$} & 
\multicolumn{3}{|c}{$k_{\theta ,n}=482$} \\ \hline
\multicolumn{1}{r|}{Test / Size} & $\mathbf{1\%}$ & $\mathbf{5\%}$ & $%
\mathbf{10\%}$ & $\mathbf{1\%}$ & $\mathbf{5\%}$ & $\mathbf{10\%}$ & $%
\mathbf{1\%}$ & $\mathbf{5\%}$ & $\mathbf{10\%}$ & $\mathbf{1\%}$ & $\mathbf{%
5\%}$ & $\mathbf{10\%}$ & $\mathbf{1\%}$ & $\mathbf{5\%}$ & $\mathbf{10\%}$
& $\mathbf{1\%}$ & $\mathbf{5\%}$ & $\mathbf{10\%}$ \\ \hline
\multicolumn{1}{l|}{$p$-Max-Test} & .010 & .042 & .116 & .013 & .053 & .103
& .008 & .055 & \multicolumn{1}{l|}{.096} & .923 & .981 & .993 & 1.00 & 1.00
& 1.00 & 1.00 & 1.00 & 1.00 \\ 
\multicolumn{1}{l|}{$p$-Max-t-Test} & .007 & .058 & .111 & .013 & .061 & .117
& .012 & .059 & \multicolumn{1}{l|}{.112} & .946 & .986 & .994 & 1.00 & 1.00
& 1.00 & 1.00 & 1.00 & 1.00 \\ 
\multicolumn{1}{l|}{$dbl$-Max-Test} & .014 & .054 & .108 & .000 & .022 & .068
& .021 & .078 & .176 & .012 & .048 & .092 & .020 & .044 & .100 & .036 & .180
& .332 \\ 
\multicolumn{1}{l|}{$dbl$-Max-t-Test} & .014 & .032 & .082 & .004 & .044 & 
.090 & .028 & .084 & .196 & .016 & .056 & .084 & .132 & .188 & .220 & .804 & 
.900 & .904 \\ 
Wald & .000 & .000 & \multicolumn{1}{l|}{.018} & \multicolumn{1}{|c}{-} & 
\multicolumn{1}{c}{-} & \multicolumn{1}{c|}{-} & \multicolumn{1}{|c}{-} & 
\multicolumn{1}{c}{-} & \multicolumn{1}{c|}{-} & .002 & .113 & .341 & 
\multicolumn{1}{|c}{-} & \multicolumn{1}{c}{-} & \multicolumn{1}{c|}{-} & 
\multicolumn{1}{|c}{-} & \multicolumn{1}{c}{-} & \multicolumn{1}{c}{-} \\ 
\hline\hline
\multicolumn{1}{c|}{} & \multicolumn{9}{|c|}{$n=250$} & \multicolumn{9}{|c}{$%
n=250$} \\ \hline
& \multicolumn{3}{|c|}{$k_{\theta ,n}=35$} & \multicolumn{3}{|c|}{$k_{\theta
,n}=1144$} & \multicolumn{3}{|c|}{$k_{\theta ,n}=1250$} & 
\multicolumn{3}{|c|}{$k_{\theta ,n}=35$} & \multicolumn{3}{|c|}{$k_{\theta
,n}=1144$} & \multicolumn{3}{|c}{$k_{\theta ,n}=1250$} \\ \hline
\multicolumn{1}{r|}{} & $\mathbf{1\%}$ & $\mathbf{5\%}$ & $\mathbf{10\%}$ & $%
\mathbf{1\%}$ & $\mathbf{5\%}$ & $\mathbf{10\%}$ & $\mathbf{1\%}$ & $\mathbf{%
5\%}$ & $\mathbf{10\%}$ & $\mathbf{1\%}$ & $\mathbf{5\%}$ & $\mathbf{10\%}$
& $\mathbf{1\%}$ & $\mathbf{5\%}$ & $\mathbf{10\%}$ & $\mathbf{1\%}$ & $%
\mathbf{5\%}$ & $\mathbf{10\%}$ \\ \hline
\multicolumn{1}{l|}{$p$-Max-Test} & .010 & .056 & .106 & .012 & .058 & .112
& .007 & .046 & \multicolumn{1}{l|}{.095} & 1.00 & 1.00 & 1.00 & 1.00 & 1.00
& 1.00 & 1.00 & 1.00 & 1.00 \\ 
\multicolumn{1}{l|}{$p$-Max-t-Test} & .012 & .055 & .106 & .015 & .061 & .116
& .009 & .057 & \multicolumn{1}{l|}{.097} & 1.00 & 1.00 & 1.00 & 1.00 & 1.00
& 1.00 & 1.00 & 1.00 & 1.00 \\ 
Wald & .000 & .015 & \multicolumn{1}{l|}{.066} & \multicolumn{1}{|c}{-} & 
\multicolumn{1}{c}{-} & \multicolumn{1}{c|}{-} & \multicolumn{1}{|c}{-} & 
\multicolumn{1}{c}{-} & \multicolumn{1}{c|}{-} & .961 & .996 & 1.00 & 
\multicolumn{1}{|c}{-} & \multicolumn{1}{c}{-} & \multicolumn{1}{c|}{-} & 
\multicolumn{1}{|c}{-} & \multicolumn{1}{c}{-} & \multicolumn{1}{c}{-} \\ 
\hline\hline
\multicolumn{1}{c|}{} & \multicolumn{9}{|c|}{$n=500$} & \multicolumn{9}{|c}{$%
n=500$} \\ \hline
& \multicolumn{3}{|c|}{$k_{\theta ,n}=35$} & \multicolumn{3}{|c|}{$k_{\theta
,n}=2381$} & \multicolumn{3}{|c|}{$k_{\theta ,n}=2500$} & 
\multicolumn{3}{|c|}{$k_{\theta ,n}=35$} & \multicolumn{3}{|c|}{$k_{\theta
,n}=2381$} & \multicolumn{3}{|c}{$k_{\theta ,n}=2500$} \\ \hline
\multicolumn{1}{r|}{Test / Size} & $\mathbf{1\%}$ & $\mathbf{5\%}$ & $%
\mathbf{10\%}$ & $\mathbf{1\%}$ & $\mathbf{5\%}$ & $\mathbf{10\%}$ & $%
\mathbf{1\%}$ & $\mathbf{5\%}$ & $\mathbf{10\%}$ & $\mathbf{1\%}$ & $\mathbf{%
5\%}$ & $\mathbf{10\%}$ & $\mathbf{1\%}$ & $\mathbf{5\%}$ & $\mathbf{10\%}$
& $\mathbf{1\%}$ & $\mathbf{5\%}$ & $\mathbf{10\%}$ \\ \hline
\multicolumn{1}{l|}{$p$-Max-Test} & .012 & 053 & .105 & .011 & .055 & .103 & 
.007 & .052 & \multicolumn{1}{l|}{.103} & 1.00 & 1.00 & 1.00 & 1.00 & 1.00 & 
1.00 & 1.00 & 1.00 & 1.00 \\ 
\multicolumn{1}{l|}{$p$-Max-t-Test} & .012 & .051 & .107 & .009 & .054 & .111
& .009 & .054 & \multicolumn{1}{l|}{.111} & 1.00 & 1.00 & 1.00 & 1.00 & 1.00
& 1.00 & 1.00 & 1.00 & 1.00 \\ 
Wald & .003 & .043 & \multicolumn{1}{l|}{.076} & \multicolumn{1}{|c}{-} & 
\multicolumn{1}{c}{-} & \multicolumn{1}{c|}{-} & \multicolumn{1}{|c}{-} & 
\multicolumn{1}{c}{-} & \multicolumn{1}{c|}{-} & .998 & 1.00 & 1.00 & 
\multicolumn{1}{|c}{-} & \multicolumn{1}{c}{-} & \multicolumn{1}{c|}{-} & 
\multicolumn{1}{|c}{-} & \multicolumn{1}{c}{-} & \multicolumn{1}{c}{-} \\ 
\hline\hline
\end{tabular}%
}
\end{center}

{\small $p$ = parsimonious; $dbl$ = de-biased Lasso. $x_{t}$ are unbounded. $k_{\delta }=0$. All test p-values are bootstrapped, based on 1,000 ($p$-max, Wald) or 250 ($dbl$-max) independently drawn samples.} 
\end{sidewaystable}

\clearpage

\begin{sidewaystable}[tbp]
\caption{Rejection Frequencies under $H_{1}(ii)$ and $H_{1}(iii)$} \label%
{table:max_wald_rejH12H13}

\begin{center}
{\small 
\begin{tabular}{c|lll|lll|llllll|lll|lll}
\hline
\  & \multicolumn{9}{|c}{$H_{1}(ii)$ : $\theta _{0,i}$ $=$ $.0002i/k_{\theta ,n}$
for $i$ $=$ $1,...,k_{\theta ,n}$} & \multicolumn{9}{|c}{$H_{1}(iii)$ : $%
\theta _{0,i}$ $=$ $.00015$ for $i$ $=$ $1,...,10$} \\ \hline
& \multicolumn{9}{|c}{$n=100$} & \multicolumn{9}{|c}{$n=100$} \\ \hline\hline
\multicolumn{1}{l|}{} & \multicolumn{3}{|c|}{$k_{\theta ,n}=35$} & 
\multicolumn{3}{|c|}{$k_{\theta ,n}=200$} & \multicolumn{3}{|c}{$k_{\theta
,n}=482$} & \multicolumn{3}{|c|}{$k_{\theta ,n}=35$} & \multicolumn{3}{|c|}{$%
k_{\theta ,n}=200$} & \multicolumn{3}{|c}{$k_{\theta ,n}=482$} \\ \hline
\multicolumn{1}{r|}{Test / Size} & $\mathbf{1\%}$ & $\mathbf{5\%}$ & $%
\mathbf{10\%}$ & $\mathbf{1\%}$ & $\mathbf{5\%}$ & $\mathbf{10\%}$ & $%
\mathbf{1\%}$ & $\mathbf{5\%}$ & $\mathbf{10\%}$ & \multicolumn{1}{|l}{$%
\mathbf{1\%}$} & $\mathbf{5\%}$ & $\mathbf{10\%}$ & $\mathbf{1\%}$ & $%
\mathbf{5\%}$ & $\mathbf{10\%}$ & $\mathbf{1\%}$ & $\mathbf{5\%}$ & $\mathbf{%
10\%}$ \\ \hline
\multicolumn{1}{l|}{$p$-Max-Test} & 1.00 & 1.00 & 1.00 & 1.00 & 1.00 & 1.00
& 1.00 & 1.00 & 1.00 & \multicolumn{1}{|l}{.934} & .988 & .993 & 1.00 & 1.00
& 1.00 & 1.00 & 1.00 & 1.00 \\ 
\multicolumn{1}{l|}{$p$-Max-t-Test} & 1.00 & 1.00 & 1.00 & 1.00 & 1.00 & 1.00
& 1.00 & 1.00 & 1.00 & \multicolumn{1}{|l}{.954} & .989 & .993 & 1.00 & 1.00
& 1.00 & 1.00 & 1.00 & 1.00 \\ 
\multicolumn{1}{l|}{$dbl$-Max-Test} & .016 & .060 & .116 & .188 & .436 & .596
& .468 & .748 & .864 & \multicolumn{1}{|l}{.016} & .048 & .092 & .000 & .036
& .084 & .020 & .072 & .120 \\ 
\multicolumn{1}{l|}{$dbl$-Max-t-Test} & .008 & .036 & .088 & .380 & .536 & 
.584 & .503 & .810 & .900 & \multicolumn{1}{|l}{.012} & .036 & .076 & .000 & 
.008 & .016 & .016 & .024 & .028 \\ 
\multicolumn{1}{l|}{Wald} & .542 & .982 & \multicolumn{1}{l|}{.998} & 
\multicolumn{1}{|c}{-} & \multicolumn{1}{c}{-} & \multicolumn{1}{c|}{-} & 
\multicolumn{1}{|c}{-} & \multicolumn{1}{c}{-} & \multicolumn{1}{c|}{-} & 
\multicolumn{1}{|l}{.003} & .101 & .307 & \multicolumn{1}{|c}{-} & 
\multicolumn{1}{c}{-} & \multicolumn{1}{c|}{-} & \multicolumn{1}{|c}{-} & 
\multicolumn{1}{c}{-} & \multicolumn{1}{c}{-} \\ \hline\hline
& \multicolumn{9}{|c}{$n=250$} & \multicolumn{9}{|c}{$n=250$} \\ \hline
\multicolumn{1}{l|}{} & \multicolumn{3}{|c|}{$k_{\theta ,n}=35$} & 
\multicolumn{3}{|c|}{$k_{\theta ,n}=1144$} & \multicolumn{3}{|c}{$k_{\theta
,n}=1250$} & \multicolumn{3}{|c|}{$k_{\theta ,n}=35$} & \multicolumn{3}{|c|}{%
$k_{\theta ,n}=1144$} & \multicolumn{3}{|c}{$k_{\theta ,n}=1250$} \\ \hline
\multicolumn{1}{r|}{Test / Size} & $\mathbf{1\%}$ & $\mathbf{5\%}$ & $%
\mathbf{10\%}$ & $\mathbf{1\%}$ & $\mathbf{5\%}$ & $\mathbf{10\%}$ & $%
\mathbf{1\%}$ & $\mathbf{5\%}$ & $\mathbf{10\%}$ & \multicolumn{1}{|l}{$%
\mathbf{1\%}$} & $\mathbf{5\%}$ & $\mathbf{10\%}$ & $\mathbf{1\%}$ & $%
\mathbf{5\%}$ & $\mathbf{10\%}$ & $\mathbf{1\%}$ & $\mathbf{5\%}$ & $\mathbf{%
10\%}$ \\ \hline
\multicolumn{1}{l|}{$p$-Max-Test} & 1.00 & 1.00 & 1.00 & 1.00 & 1.00 & 1.00
& 1.00 & 1.00 & 1.00 & \multicolumn{1}{|l}{.998} & 1.00 & 1.00 & 1.00 & 1.00
& 1.00 & 1.00 & 1.00 & 1.00 \\ 
\multicolumn{1}{l|}{$p$-Max-t-Test} & 1.00 & 1.00 & 1.00 & 1.00 & 1.00 & 1.00
& 1.00 & 1.00 & 1.00 & \multicolumn{1}{|l}{1.00} & 1.00 & 1.00 & 1.00 & 1.00
& 1.00 & 1.00 & 1.00 & 1.00 \\ 
\multicolumn{1}{l|}{Wald} & 1.00 & 1.00 & 1.00 & \multicolumn{1}{|c}{-} & 
\multicolumn{1}{c}{-} & \multicolumn{1}{c|}{-} & \multicolumn{1}{|c}{-} & 
\multicolumn{1}{c}{-} & \multicolumn{1}{c|}{-} & \multicolumn{1}{|l}{.730} & 
.944 & .980 & \multicolumn{1}{|c}{-} & \multicolumn{1}{c}{-} & 
\multicolumn{1}{c|}{-} & \multicolumn{1}{|c}{-} & \multicolumn{1}{c}{-} & 
\multicolumn{1}{c}{-} \\ \hline\hline
& \multicolumn{9}{|c}{$n=500$} & \multicolumn{9}{|c}{$n=500$} \\ \hline
\multicolumn{1}{l|}{} & \multicolumn{3}{|c|}{$k_{\theta ,n}=35$} & 
\multicolumn{3}{|c|}{$k_{\theta ,n}=2381$} & \multicolumn{3}{|c}{$k_{\theta
,n}=2500$} & \multicolumn{3}{|c|}{$k_{\theta ,n}=35$} & \multicolumn{3}{|c|}{%
$k_{\theta ,n}=2381$} & \multicolumn{3}{|c}{$k_{\theta ,n}=2500$} \\ \hline
\multicolumn{1}{r|}{Test / Size} & $\mathbf{1\%}$ & $\mathbf{5\%}$ & $%
\mathbf{10\%}$ & $\mathbf{1\%}$ & $\mathbf{5\%}$ & $\mathbf{10\%}$ & $%
\mathbf{1\%}$ & $\mathbf{5\%}$ & $\mathbf{10\%}$ & \multicolumn{1}{|l}{$%
\mathbf{1\%}$} & $\mathbf{5\%}$ & $\mathbf{10\%}$ & $\mathbf{1\%}$ & $%
\mathbf{5\%}$ & $\mathbf{10\%}$ & $\mathbf{1\%}$ & $\mathbf{5\%}$ & $\mathbf{%
10\%}$ \\ \hline
\multicolumn{1}{l|}{$p$-Max-Test} & 1.00 & 1.00 & 1.00 & 1.00 & 1.00 & 1.00
& 1.00 & 1.00 & 1.00 & \multicolumn{1}{|l}{1.00} & 1.00 & 1.00 & 1.00 & 1.00
& 1.00 & 1.00 & 1.00 & 1.00 \\ 
\multicolumn{1}{l|}{$p$-Max-t-Test} & 1.00 & 1.00 & 1.00 & 1.00 & 1.00 & 1.00
& 1.00 & 1.00 & 1.00 & \multicolumn{1}{|l}{1.00} & 1.00 & 1.00 & 1.00 & 1.00
& 1.00 & 1.00 & 1.00 & 1.00 \\ 
\multicolumn{1}{l|}{Wald} & 1.00 & 1.00 & 1.00 & \multicolumn{1}{|c}{-} & 
\multicolumn{1}{c}{-} & \multicolumn{1}{c|}{-} & \multicolumn{1}{|c}{-} & 
\multicolumn{1}{c}{-} & \multicolumn{1}{c|}{-} & 1.00 & 1.00 & 1.00 & 
\multicolumn{1}{|c}{-} & \multicolumn{1}{c}{-} & \multicolumn{1}{c|}{-} & 
\multicolumn{1}{|c}{-} & \multicolumn{1}{c}{-} & \multicolumn{1}{c}{-} \\ 
\hline\hline
\end{tabular}
}
\end{center}

{\small $p$ = parsimonious; $dbl$ = de-biased Lasso. $x_{t}$ are unbounded. $k_{\delta }=0$. All test p-values are bootstrapped, based on 1,000 ($p$-max, Wald) or 250 ($dbl$-max) independently drawn samples.} 
\end{sidewaystable}

\end{document}